\pdfoutput=1 
\documentclass[leqno]{article}
\usepackage[french,english]{babel}
\usepackage[normal]{phaine_style}
\usepackage[margin=1.5in]{geometry}

\numberwithin{equation}{subsection}

\hypersetup{bookmarksdepth=2}

\DeclareMathOperator{\Link}{Link}

\theoremstyle{definition}

\newtheorem{comparison}[equation]{Comparison}

\crefname{comparison}{Comparison}{Comparisons}
\Crefname{comparison}{Comparison}{Comparisons}

\newcommand{\StrTop}{\categ{StrTop}^{\ex}}
\newcommand{\Str}{\categ{Str}}
\newcommand{\Joyal}{\mathrm{Joy}}
\newcommand{\sSetJoyal}{\sSet^{\,\Joyal}}
\newcommand{\TopP}{\Top_{/P}}
\newcommand{\ex}{\mathrm{ex}}
\newcommand{\seg}{\mathrm{seg}}
\newcommand{\TopPex}{\Top_{/P}^{\ex}}
\newcommand{\TopPseg}{\Top_{/P}^{\seg}}

\newcommand{\TopB}{\Top_{/B}}

\newcommand{\sSetP}{\sSet_{/P}}

\newcommand{\StrP}{\Str_P}
\renewcommand{\Spc}{\Gpd_{\infty}}
\newcommand{\sdop}[1]{\sd(#1)^{\op}}
\newcommand{\sSetPproj}{\Fun(\sdop{P},\sSet)^{\proj}}
\newcommand{\sSetDeltaPproj}{\Fun(\DDelta_{/P}^{\op},\sSet)^{\proj}}
\DeclareMathOperator{\Pair}{Pair}

\newcommand{\Deccategory}{\categ{Déc}}
\newcommand{\DecP}{\Deccategory_{P}}
\newcommand{\Con}{\categ{Con}}

\renewcommand{\W}{\mathrm{W}}

\newcommand{\nv}{\mathrm{nv}}
\newcommand{\IH}{\mathrm{IH}}
\newcommand{\LH}{\mathrm{LH}}
\newcommand{\Inv}{\IH^{\nv}}

\newcommand{\finverse}{f^{-1}}

\newcommand{\NerveP}{\Nerve_{P}}
\renewcommand{\D}{\mathrm{D}}

\newcommand{\LP}{\mathrm{L}_P}

\newcommand{\DP}{\D_P}
\newcommand{\EP}{\mathrm{E}_P}
\newcommand{\JP}{\mathrm{J}_P}

\newcommand{\WP}{\mathrm{W}_P}

\newcommand{\spine}[1]{\operatorname{Spn}^{#1}}
\newcommand{\MapP}{\Map_{\kern-0.05em/P}}
\newcommand{\Deltadot}{\Delta^{\scriptscriptstyle \bullet}}
\renewcommand{\id}[1]{\operatorname{id}_{#1}}

\DeclareMathOperator{\VEx}{VEx}
\newcommand{\ExPhat}{\VEx_P}
\newcommand{\SingP}{\Sing_{P}}
\newcommand{\real}{\abs}
\newcommand{\realP}[1]{\real{#1}_{P}}

\newcommand{\ExitP}{\Exit_P}

\def\reflectedop#1#2{\mathop{\reflectbox{$#1{#2}$}}}
\renewcommand{\tensor}{\rtimes}
\renewcommand{\tensorhat}{\mathop{\widehat{\tensor}}} 
\renewcommand{\crosshat}{\mathop{\widehat{\cross}}} 

\usepackage{todonotes}
\usepackage{xargs}  
\newcommandx{\PHtodo}[2][1=]{\todo[linecolor=blue,backgroundcolor=blue!25,bordercolor=blue,#1]{#2}}

\title{On the homotopy theory of stratified spaces}

\author{Peter J. Haine}
\date{\today}

\begin{document}

\maketitle

\begin{abstract}
	Let $ P $ be a poset.
	We define a new homotopy theory of suitably nice $ P $-stratified topological spaces with equivalences on strata and links inverted.
	We show that the exit-path construction of MacPherson, Treumann, and Lurie defines an equivalence from our homotopy theory of $ P $-stratified topological spaces to the \category of \categories with a conservative functor to $ P $.
	This proves a stratified form of Grothendieck's homotopy hypothesis, verifying a conjecture of Ayala--Francis--Rozenblyum.
	Our homotopy theory of stratified spaces has the added benefit of capturing all examples of geometric interest: conically stratified spaces fit into our theory, and the Ayala--Francis--Tanaka--Rozenblyum homotopy theory of conically smooth stratified spaces embeds into ours.
\end{abstract}


\tableofcontents


\setcounter{section}{-1}

\section{Introduction}

Let $ T $ be a topological space.
Classically, there are two approaches to understanding the algebraic topology of $ T $.
On the one hand, we can work directly with the genuine topological space $ T $ and think about homotopy classes of maps to $ T $, homotopy groups, etc.
On the other hand, we can work with a more combinatorial object: the \textit{singular simplicial set} or \textit{fundamental \groupoid} of $ T $.
This combinatorial object parametrizes points of $ T $, paths in $ T $, homotopies of paths, etc.
Grothendieck's celebrated \textit{homotopy hypothesis}, established by Kan and Quillen \cites{MR90047}{MR0096210}{MR111032}{MR111033}{MR111036}{MR0223432}, asserts that these two perspectives are equivalent.
The equivalence between the homotopy theories of topological spaces and simplicial sets is at the foundation of modern homotopy theory for a good reason: both perspectives are of great utility. 
Being able to think about homotopy types as topological spaces gives access to many examples and structures coming from the geometry of manifolds.
On the other hand, the combinatorial framework of simplicial sets provides a rich context in which to do algebra in the setting of homotopy theory (e.g., study structured ring spectra).

The purpose of this article is to put the homotopy theory of \textit{stratified} spaces on a similarly good footing.
One peculiarity is that for a long time only a robust `combinatorial' approach to the homotopy theory of stratified spaces existed.
The combinatorial approach originates from ideas of MacPherson that came out of studying manifolds with singularities and constructible sheaves.
MacPherson observed that a constructible sheaf of vector spaces on a (nice) stratified space $ T $ is equivalent to a a representation of the \textit{exit-path} category of $ T $.
To explain this, consider the simple example of the unit interval $ [0,1] $ with stratification $ \{0\} \subset [0,1] $.
Let $ F $ be a sheaf on $ [0,1] $ that is locally constant on $ (0,1] $, i.e., constructible with respect to the stratification.
Since $ \restrict{F}{(0,1]} $ is locally constant, each path between points $ s,t \in (0,1] $ defines an isomorphism between stalks $ F_s \isomorphic F_t $. 
Using these isomorphisms and the fact that every open neighborhood of $ 0 $ intersects $ (0,1] $, one can define a \textit{specalization map} $ \fromto{F_0}{F_1} $ relating the stalks of $ F $ at $ 0 $ and $ 1 $.
These three pieces of data completely determine the constructible sheaf $ F $. 
That is, constructible sheaves on $ \{0\} \subset [0,1] $ are representations of the $ \Aup_2 $-quiver $ \bullet \to \bullet $.

More generally, if $ T $ is a topological space with a suitably nice stratification by a poset $ P $, then we can associate to $ T $ its \textit{exit-path \category} with objects points of $ T $, with $ 1 $-morphisms \textit{exit paths} flowing from lower to higher strata (and once they exit a stratum are not allowed to return), with $ 2 $-morphisms homotopies of exit-paths respecting stratifications, etc. 
The adjective `suitably nice' is quite important here because, while the construction of the fundamental \groupoid makes sense for \textit{any} topological space, if the stratification is not sufficiently nice, then exit paths can fail to suitably compose and this informal description cannot be made to actually define an \category.
This is part of an overarching problem: there does not yet exist a homotopy theory of stratified spaces that is simple to define, encapsulates examples from topology,
and has excellent formal properties.
One of the goals of this paper is to resolve this matter.

Treumann \cite{MR2575092}, Woolf \cite{MR2591969}, Lurie \cite[\HAapp{A}]{HA}, and Ayala--Francis--Rozen\-blyum \cite[\S1]{MR3941460} have all worked to realize MacPherson's exit-path construction under a variety of point-set topological assumptions.
In various forms, Ayala--Francis--Rozenblyum \cite[Conjectures 0.0.4 \& 0.0.8]{MR3941460}, Barwick, and Woolf have all conjectured that the exit-path construction defines an equivalence of \categories from suitably nice stratified spaces (with stratified homotopies inverted) to \categories with a conservative functor to a poset.
The main goal of this paper is to prove this conjecture.

There is already strong evidence for the power of having both topological and combinatorial approaches to stratified homotopy theory on the same footing.
On the one hand, Ayala--Francis--Tanaka--Rozenblyum have made extensive use of explicit topological methods to study manifolds with singularities and their factorization homology \cites{MR3916983}{MR3941460}{MR3595895}{MR3590534}.
On the other hand, in joint work with Barwick and Glasman \cite{exodromy}, we used the combinatorial perspective on stratified spaces to associate to each variety an `exit-path category' for the étale topology.
The combinatorial perspective is necessary here because `differential-topological' constructions of exit-path \categories are not available in positive-characteristic algebraic geometry.
This étale exit-path category is a powerful globalization of the absolute Galois groups of a field.
For example, it gives rise to a new and concrete definition of the étale homotopy type of Artin--Mazur--Friedlander \cites{MR0245577}{MR676809}.
Moreover, much like how the Neukirch--Uchida theorem shows that the absolute Galois group is a complete invariant of number fields \cites{MR0244211}{MR0258804}{MR0432593}, the étale exit-path category is a complete invariant of the varieties that appear in Grothendieck's anabelian conjecture \cite{MR1483108}. 


\subsection{The stratified homotopy hypothesis}

In order to state our stratified homotopy hypothesis, it is useful to get a better understanding of what kind of data should determine a stratified space.
A stratification of a topological space $ T $ by the poset $ \{0 < 1\} $ is given by a closed subspace $ T_0 \subset T $ and its open complement $ T_1 \colonequals T \sminus T_0 $.
It is natural to ask: given the topological spaces $ T_0 $ and $ T_1 $, what extra data do we need to reconstruct the stratified space $ T $ (up to stratified homotopy equivalence)?
Roughly, the answer is \textit{gluing} information called the \textit{(homotopy) link} between the $ 0 $-th and $ 1 $-st strata.
Introduced by Quinn \cite{MR928266}, the link $ \Link(T_0,T_1) $ is defined as the space of paths $ \gamma \colon \fromto{[0,1]}{T} $ such that $ \gamma(0) \in T_0 $ and for all $ s > 0 $ we have $ \gamma(s) \in T_1 $.
Evaluation at $ 0 $ and $ 1 $ define maps $ \fromto{\Link(T_0,T_1)}{T_0} $ and $ \fromto{\Link(T_0,T_1)}{T_1} $, respectively.
In nice situations, the stratified space $ T $ can be recovered as a homotopy pushout of the resulting span
\begin{equation*}
	T_0 \ot \Link(T_0,T_1) \to T_1 
\end{equation*}
\cite[\S2]{MR928266}.
For more general stratifications, the idea is that stratified spaces can be reconstructed from the data of all of their strata and all possible links relating strata. 

To state our result, we introduce some notation.
Given a poset $ P $, we write $ \TopP $ for the category of $ P $-stratified topological spaces (see \Cref{def:P-stratified_spaces}).
Lurie's exit-path construction defines a functor \smash{$ \SingP \colon \fromto{\TopP}{\sSetP} $} from $ P $-stratified spaces the category of simplicial sets over (the nerve of) $ P $ (see \cref{sec:recstrattopspace}).
Write \smash{$ \TopPex \subset \TopP $} for the full subcategory of those $ P $-stratified topological spaces $ T $ for which Lurie's exit-path simplicial set $ \SingP(T) $ is \acategory.
The \categories that arise in this way have a special property: the fibers of the structure morphism to $ P $ are given by the fundamental \groupoids of the strata.
In particular, every morphism in each fiber is invertible. 
Equivalently, the structure morphism is a conservative functor.
We write $ \StrP $ for the \category of \categories over (the nerve of) $ P $ such that the structure morphism is conservative.
We refer to $ \StrP $ as the \category of \textit{abstract $ P $-stratified homotopy types}.

Our main result provides an affirmative answer to a conjecture of Ayala--Francis--Rozenblyum \cite[Conjectures 0.0.4 \& 0.0.8]{MR3941460}:

\begin{theorem}[(stratified homotopy hypothesis; see \Cref{prop:SingPequiv,cor:globalSHH})]\label{introprop:SingPequiv}
	Let $ P $ be a poset.
	Lurie's exit-path construction defines an equivalence of \categories 
	\begin{equation*}
		\SingP \colon \TopPex\brackets{\left(\substack{\textup{maps inducing weak homotopy} \\ \textup{equivalences on strata and links}} \right)^{-1}} \equivalence \StrP 
	\end{equation*}
	from the \category obtained from $ \TopPex $ by inverting maps that induce weak homotopy equivalences on strata and links to the \category $ \StrP $ of abstract $ P $-stratified homotopy types.
\end{theorem}

\noindent In particular, every \category with a conservative functor to a poset arises as the exit-path \category of a stratified topological space.

A key consequence of our stratified homotopy hypothesis is that the Ayala--Francis--Tanaka--Rozenblyum homotopy theory of \textit{conically smooth stratified spaces} \cites{MR3941460}[\S3]{MR3590534} embeds into ours.
See \Cref{comp:AFR} for details.
A major benefit of our homotopy theory of stratified spaces is that all \textit{conically stratified topological spaces} fit into our framework (see \cref{nul:conicstratisqcat}).
\textit{Topologically stratified spaces} in the sense of Goresky--MacPherson \cite[\S1.1]{MR696691}, in particular all Whitney-stratified spaces \cites{MR2958928}{MR239613}, are conically stratified.
Thus our homotopy theory captures most, if not all, examples of differential-topological interest.
On the other hand, for a long time it was unknown whether or not every Whitney-stratified space admits a conically smooth structure \cite[Conjecture 0.0.7]{MR3941460}.
Only recently, Nocera and Volpe have shown that Whitney-stratified spaces admit conically smooth atlases \cite[Theorem 2.7]{arXiv:2105.09243}.


\subsection{Proof outline}

Proving \Cref{introprop:SingPequiv} requires working with explicit model categories of stratified topological spaces and simplicial sets and reinterpreting the \category $ \StrP $ in a way that relates to the links of stratified topological spaces.
To do this, we make use of Douteau's recent thesis and subsequent work \cites{arXiv:1908.01366}{MR4196384}{arXiv:2102.04876}{MR4224748}.
Douteau's works realizes the following key insight of Henriques \cites[\S4.7]{Henriques:thesis}{Henriques:stratifiedmodel}: while checking equivalences after passing to exit-path \categories is only reasonable for suitably nice $ P $-stratified topological spaces, a morphism of suitably nice $ P $-stratified topological spaces is an equivalence on exit-path \categories if and only if it induces an equivalence on all spaces of sections over geometric realizations of linearly ordered finite subsets of $ P $.
Moreover, the latter definition in terms of these `higher order links' works well for \textit{all} stratified topological spaces.

Precisely, let $ \sd(P) $ denote the \textit{subdivision} of $ P $, that is, the poset of linearly ordered finite subsets $ \Sigma \subset P $.
There is a right adjoint `nerve' functor 
\begin{align*}
	\Nerve_P \colon \TopP &\to \Fun(\sdop{P},\sSet) \\
	T &\mapsto [\Sigma \mapsto \Sing \MapP(\real{\Sigma},T)] 
\end{align*}
from the category of $ P $-stratified topological spaces to the category of simplicial presheaves on $ \sd(P) $; see \Cref{ntn:NPandLP}.
Douteau proves that the projective model structure transfers to $ \TopP $ along $ \Nerve_P $, so that a morphism $ f \colon \fromto{T}{S} $ in $ \TopP $ is a weak equivalence (resp., fibration) if and only if for every $ \Sigma \in \sd(P) $, the induced map
\begin{equation*}
	\fromto{\MapP(\real{\Sigma},T)}{\MapP(\real{\Sigma},S)}
\end{equation*}
on topological spaces of sections over the realization of $ \Sigma $ is a weak homotopy equivalence (resp., Serre fibration).
Even better, the resulting Quillen adjunction \smash{$ \Fun(\sdop{P},\sSet) \rightleftarrows \TopP $} is a simplicial Quillen equivalence of combinatorial simplicial model categories (\Cref{thm:Douteaumain}).
This means that the underlying \category of the Douteau--Henriques model structure on $ \TopP $ is equivalent to the \category \smash{$ \Fun(\sdop{P},\Spc) $} of presheaves of \groupoids on the subdivision $ \sd(P) $. 

In joint work with Barwick and Glasman, we proved that a similarly-defined `nerve' functor
\begin{align*}
	\StrP &\to \Fun(\sdop{P},\Spc) \\
	X &\mapsto [\Sigma \mapsto \Fun_{/P}(\Sigma,X)] 
\end{align*}
expresses the \category $ \StrP $ of abstract $ P $-stratified homotopy types as an accessible localization of $ \Fun(\sdop{P},\Spc) $ \cite[Theorem 2.7.4]{exodromy}.
We also explicitly identify the essential image of this functor (see \cref{subsec:decollage}). 
Combining these two works we show:

\begin{theorem}[(\Cref{thm:TopPpresentsStrP})]\label{mainthm:TopPpresentsStrP}
	Let $ P $ be a poset.
	Then the \category $ \StrP $ of abstract $ P $-stratified homotopy types is equivalent to an $ \upomega $-accessible localization of the underlying \category of the combinatorial simplicial model category $ \TopP $ in the Douteau--Henriques model structure.
\end{theorem}

\noindent That is, the \category $ \StrP $ can be obtained from the ordinary category \smash{$ \TopP $} of $ P $-stratified topological spaces by inverting a class of weak equivalences.

To prove \Cref{introprop:SingPequiv}, we need to write down functors from homotopy theories of stratified topological spaces into $ \StrP $; to do so, it is convenient to present $ \StrP $ as the underlying \category of a model category.
Consider the left Bousfield localization of the Joyal model structure inherited on the overcategory \smash{$ \sSetP $} obtained by inverting all simplicial homotopies \smash{$ \fromto{X \cross \Delta^1}{Y} $} respecting the stratifications of $ X $ and $ Y $ by $ P $ (see \Cref{def:Joyal-Kan}).
We call the resulting model structure on $ \sSetP $ the \textit{Joyal--Kan} model structure. 
Surprisingly, the Joyal--Kan model structure shares many of the excellent formal properties of the Kan model structure that the Joyal model structure lacks; namely, the Joyal--Kan model structure is \textit{simplicial}.
The following results summarize the main features of the Joyal--Kan model structure and its relation to the \category $ \StrP $.

\begin{theorem}\label{thm:modstruct}
	Let $ P $ be a poset.
	\begin{enumerate}[label=\stlabel{thm:modstruct}]
		\item The Joyal--Kan model structure on $ \sSetP $ is left proper, combinatorial, and simplicial.
		The cofibrations are the monomorphisms and the fibrant objects are the \categories $ X $ over $ P $ such that the structure morphism $ \fromto{X}{P} $ is a conservative functor (\Cref{prop:JKexists,prop:equivfibrant,thm:JKissimplicial}).

		\item If $ f \colon \fromto{X}{Y} $ a morphism in $ \sSetP $ and all of the fibers of $ X $ and $ Y $ over points of $ P $ are Kan complexes (e.g., $ X $ and $ Y $ are fibrant objects), then $ f $ is an equivalence in the Joyal--Kan model structure if and only if $ f $ is an equivalence in the Joyal model structure (\Cref{prop:WEfiberwiseKan}).
	\end{enumerate} 
\end{theorem} 

\begin{theorem}[(\Cref{cor:underlying})]\label{thm:underlying}
	Let $ P $ be a poset.
	Then the underlying \category of the Joyal--Kan model structure on $ \sSetP $ is the \category $ \StrP $ of abstract $ P $-stratified homotopy types.
\end{theorem} 

The Joyal--Kan model structure gives us the explicit control we need to prove the straitifed homotopy hypothesis (\Cref{introprop:SingPequiv}).


\begin{acknowledgments}
	We would like to thank Clark Barwick for countless insights and conversations about the material in this text.
	We would also like to thank Alex Sear for her generous hospitality during his visits to the University of Edinburgh, during which much of this work took place.

	For a long time we had a sketch of an argument proving the main theorem, and reading Sylvain Douteau's thesis provided the key to complete it.
	We thank Sylvain most heartily for generously sharing his work with us and discussing the results of his thesis.
	We also thank David Chataur for his insightful correspondence about the homotopy theory of stratified spaces.

	We are grateful to Stephen Nand-Lal and Jon Woolf for explaining their work on the homotopy theory of stratified topological spaces, and in particular to Stephen for sharing his thesis with us.

	We gratefully acknowledge support from the MIT Dean of Science Fellowship, the NSF Graduate Research Fellowship under Grant \#112237, UC President's Postdoctoral Fellowship, and NSF Mathematical Sciences Postdoctoral Research Fellowship under Grant \#DMS-2102957. 
\end{acknowledgments}


\subsection{Terminology \& notations}

\begin{nul}\label{nul:terminology}
	We use the language and tools of higher category theory, particularly in the model of \textit{quasicategories}, as defined by Boardman--Vogt and developed by Joyal and Lurie \cites{HTT}{HA}.
	\begin{enumerate}
		\item We write $ \sSet $ for the category of simplicial sets and $ \Map \colon \fromto{\sSet^{\op} \cross \sSet}{\sSet} $ for the internal-$ \Hom $ in simplicial sets.

		\item To avoid confusion, we call weak equivalences in the Joyal model structure on $ \sSet $ \defn{Joyal equivalences} and we call weak equivalences in the Kan model structure on $ \sSet $ \defn{Kan equivalences}.

		\item We write $ \sSetJoyal $ for the model category of simplicial sets in the Joyal model structure.

		\item An \textit{\category} here will always mean \textit{quasicategory}; we write $ \Cat_{\infty} $ for the \category of \categories.
		We write $ \Spc \subset \Cat_{\infty} $ for the \category of \groupoids, i.e., the \textit{\category of spaces}.
		In order not to overload the term `space', we use `\groupoid' to refer to homotopy types, and `space' only in reference to topological spaces.

		\item If $ C $ is an ordinary category, we simply write $ C \in \sSet $ for its nerve.

		\item For an \category $ C $, we write $ C^{\equivalent} \subset C $ for the maximal sub-\groupoid contained in $ C $.

		\item Let $ C $ be an \category and $ W \subset \Mor(C) $ a class of morphisms in $ C $.
		We write $ C[W^{-1}] $ for the \defn{localization of $ C $ at $ W $}, i.e., the initial \category equipped with a functor $ \fromto{C}{C[W^{-1}]} $ that sends morphisms in $ W $ to equivalences \cite[\S7.1]{MR3931682}.

		\item The \defn{underlying quasicategory} of a simplicial model category $ \AA $ is the homotopy-coherent nerve of the full subcategory spanned by the fibrant--cofibrant objects (which forms a fibrant simplicial category).
		In the quasicategory model, this is a presentation of the localization of $ \AA $ at its class of weak equivalences.

		\item For every integer $ n \geq 0 $, we write $ [n] $ for the linearly ordered poset $ \{0 < 1 < \cdots < n\} $ of cardinality $ n + 1 $ (whose nerve is the simplicial set $ \Delta^n $).

		\item We denote an adjunction of categories or \categories by $ \adjto{F}{C}{D}{G} $, where $ F $ is the left adjoint and $ G $ is the right adjoint.

		\item To fix a convenient category of topological spaces, we write $ \Top $ for the category of \textit{numerically generated} topological spaces (also called \textit{$ \Delta $-generated} or \textit{$ I $-generated} topological spaces) \cites{DuggerDelta}{MR2425555}{HaraguchiThesis}[\S3]{MR3289294}{MR3884529}, and use the term `topological space' to mean `numerically generated topological space'. 
		For the present work, the category of numerically generated topological spaces is preferable to the more standard category of compactly generated weakly Hausdorff topological spaces \cite[Chapter 5]{MR1702278}: a poset in the Alexandroff topology is weakly Hausdorff if and only if it is discrete, whereas every poset in is numerically generated.
	\end{enumerate}
\end{nul}

\begin{definition}
	Let $ P $ be a poset.
	The category of \defn{$ P $-stratified simplicial sets} is the overcategory $ \sSetP $ of simplicial sets over (the nerve of) $ P $.
	Given a $ P $-stratified simplicial set $ f \colon \fromto{X}{P} $ and point $ p \in P$, we write \smash{$ X_p \colonequals \finverse(p) $} for the \smash{\defn{$ p $-th stratum}} of $ X $.
\end{definition}

\begin{notation}
	Let $ P $ be a poset.
	Write $ - \tensor_P - \colon \fromto{\sSetP \cross \sSet}{\sSetP} $ for the standard tensoring of $ \sSetP $ over $ \sSet $, defined on objects by sending an object $ X \in \sSetP $ and a simplicial set $ K \in \sSet $ to the product $ X \tensor_P K \colonequals X \cross K $ in $ \sSet $ with structure morphism induced by the projection $ \fromto{X \cross K}{X} $.
    When unambiguous we write $ \tensor $ rather than $ \tensor_P $, leaving the poset $ P $ implicit. 

    We write \smash{$ \MapP \colon \fromto{\sSetP^{\op} \cross \sSetP}{\sSet} $} for the standard simplicial enrichment, whose assignment on objects is given by
    \begin{equation*}
    	\MapP(X,Y) \colonequals \sSetP(X \tensor_P \Deltadot,Y) \comma
    \end{equation*}
    and the assignment on morphisms is the obvious one.
\end{notation}


\section{Abstract stratified homotopy types, décollages, \& stratified topological spaces}\label{sec:thmA}

This section is dedicated to the proof of \Cref{mainthm:TopPpresentsStrP}.


\subsection{Abstract stratified homotopy types as décollages}\label{subsec:decollage}

In work with Barwick and Glasman \cite[\S\S2.6 \& 2.7]{exodromy}, we gave a complete Segal space style description of the \category of abstract $ P $-stratified homotopy types. 
In this subsection we recall this work and, for completeness, include a proof of the main comparison result (\Cref{thm:nerveequiv}).

\begin{definition}
	Let $ P $ be a poset.
	The \category $ \StrP $ of \defn{abstract $ P $-stratified homotopy types} is the full subcategory of the overcategory $ \Cat_{\infty,/P} $ spanned by those \categories over $ P $ with conservative structure morphism $ \fromto{C}{P} $.

	Note that the mapping \groupoid \smash{$ \Map_{\StrP}(X,Y) $} coincides with the \category \smash{$ \Fun_{/P}(X,Y) $} of functors $ \fromto{X}{Y} $ over $ P $.
\end{definition}

\begin{recollection}\label{rec:Segal_spaces}
	\Acategory can be modeled as a simplicial \groupoid.
	There is a nerve functor $ \Nerve \colon \fromto{\Cat_{\infty}}{\Fun(\Deltaop,\Spc)} $ defined by
	\begin{equation*}
		\Nerve(C)_m \coloneq \Fun(\Delta^m,C)^{\equivalent} \period
	\end{equation*}
	The simplicial \groupoid $ \Nerve(C) $ is an example of what Rezk called a \defn{complete Segal space} \cite{MR1804411}, i.e., a functor $ F \colon \fromto{\Deltaop}{\Spc} $ satisfying the following conditions:
	\begin{enumerate}[label=\stlabel{rec:Segal_spaces}]
		\item \textit{Segal condition:} For each integer $ m \geq 1 $, the natural map
		\begin{equation*}
			F_m \to F\{0 < 1\} \crosslimits_{F\{1\}} F\{1 < 2\} \crosslimits_{F\{2\}} \cdots \crosslimits_{F\{m-1\}} F\{m-1 < m\}
		\end{equation*}
		is an equivalence.

		\item \textit{Completeness condition:} The natural morphism
		\begin{equation*}
			F_0 \to F_3 \crosslimits_{F\{0 < 2\} \cross F\{1 < 3\}} (F_0 \cross F_0) 
		\end{equation*}
		is an equivalence.
	\end{enumerate}

	Joyal and Tierney showed that the nerve is fully faithful with essential image the full subcategory $\CSS$ spanned by the complete Segal spaces \cite{MR2342834}.
\end{recollection}

We now give an analogous description of $ \StrP $.

\begin{notation}
	Let $ P $ be a poset.
	We write $ \sd(P) $ for the poset of nonempty linearly ordered finite subsets $ \Sigma \subset P $ ordered by containment.
	We call $ \sd(P) $ the \defn{subdivision} of $ P $.
	We call a nonempty linearly ordered finite subset $ \Sigma \subset P $ of $ P $ a \defn{chain} in $ P $.
\end{notation}

\begin{definition}
	Let $ P $ be a poset.
	A functor $ F \colon \fromto{\sdop{P}}{\Spc} $ is a \defn{décollage} (over $P$) if and only if, for every chain $ \{p_0 < \cdots < p_m\} \subset P$, the natural map
	\begin{equation*}
		F\{p_0 < \cdots < p_m\} \to F\{p_0 < p_1\} \crosslimits_{F\{p_1\}} F\{p_1 < p_2\} \crosslimits_{F\{p_2\}} \cdots \crosslimits_{F\{p_{m-1}\}} F\{p_{m-1} < p_m\}
	\end{equation*}
	is an equivalence.
	We write
	\begin{equation*}
		\DecP \subset \Fun(\sdop{P},\Spc)
	\end{equation*}
	for the full subcategory spanned by the décollages.
	Note that $ \DecP $ is closed under limits and filtered colimits in $ \Fun(\sdop{P},\Spc) $.
\end{definition}

\noindent A nerve style construction provides an equivalence $ \equivto{\StrP}{\DecP} $.

\begin{construction}\label{cnstr:nerveofstratifiedspace} 
	Let $P$ be a poset. 
	We have a fully faithful functor $ \incto{\sd(P)}{\StrP} $ given by regarding a chain $ \Sigma $ as an \category over $ P $ via the inclusion $ \incto{\Sigma}{P} $. 
	Define a functor
	\begin{align*}
		\Nerve_P \colon \StrP &\to \Fun(\sdop{P},\Spc) \\
	\intertext{by the assignment}
		X &\mapsto [\goesto{\Sigma}{\Map_{\StrP}(\Sigma,X)}] \period
	\end{align*}
\end{construction}

\begin{nul}
	Since a chain $ \{p_0 < \cdots < p_n\} $ is the iterated pushout
	\begin{equation*}
		\{p_0 < p_1\} \union^{\{p_1\}} \cdots \union^{\{p_{n-1}\}} \{p_{n-1} < p_n\}
	\end{equation*}
	in $ \StrP $, the functor $ \Nerve_P $ lands in the full subcategory $ \DecP $.
\end{nul}

\begin{theorem}[{\cite[Theorem 2.7.4]{exodromy}}]\label{thm:nerveequiv}
	For any poset $ P $, the functor $ \Nerve_P \colon \fromto{\StrP}{\DecP} $ is an equivalence of \categories.
\end{theorem}

\begin{proof}
	Let $ \DDelta_{/P} $ denote the category of simplices of $ P $.
	The Joyal--Tierney Theorem \cite{MR2342834} implies that the nerve functor
	\begin{align*}
		\Cat_{\infty,/P} &\to \Fun(\Deltaop,\Spc)_{/\Nerve(P)}\simeq\Fun(\DDelta_{/P}^{\op},\Spc) \\ 
		X &\mapsto [\Sigma \mapsto \Fun_{/P}(\Sigma,X)^{\equivalent}]
	\end{align*}
	is fully faithful, with essential image \smash{$ \CSS_{/\Nerve(P)} $} those functors \smash{$ \fromto{\DDelta_{/P}^{\op}}{\Spc} $} that satisfy both the Segal condition and the completeness condition.
	Now notice that left Kan extension along the inclusion $ \incto{\sd(P)}{\DDelta_{/P}}$ defines a fully faithful functor $\incto{\DecP}{\CSS_{/\Nerve(P)}}$ whose essential image consists of those complete Segal spaces $\fromto{C}{\Nerve(P)}$ such that for any $ p\in P $, the complete Segal space $ C_p $ is an \groupoid.
\end{proof}

\begin{nul}\label{nul:StrPisanaccessibleloc}
	Since $ \StrP $ is presentable and $ \DecP \subset \Fun(\sdop{P},\Spc) $ is closed under limits and filtered colimits, the Adjoint Functor Theorem shows that the nerve expresses the \category $ \StrP $ as an $ \upomega $-accessible localization of $ \Fun(\sdop{P},\Spc) $.
\end{nul}

\begin{nul}\label{nul:equivcheckedonstratandlink}
	\Cref{thm:nerveequiv} implies that equivalences in $ \StrP $ are checked on strata and \textit{links}. 
	That is, a morphism $ f \colon \fromto{X}{Y} $ in $ \StrP $ is an equivalence if and only if $ f $ induces an equivalence on strata and for each pair $ p < q $ in $ P $, the induced map on \textit{links}
	\begin{equation*}
		\fromto{\Map_{\StrP}(\{p < q\}, X)}{\Map_{\StrP}(\{p < q\}, Y)}
	\end{equation*}
	is an equivalence in $ \Spc $.
	(This can also be proven directly without appealing to \Cref{thm:nerveequiv}.)
\end{nul}


\subsection{Recollections on stratified topological spaces}\label{sec:recstrattopspace}

We now recall the relationship between $ P $-stratified topological spaces and $ P $-stratified simplicial sets.
Recall that we write $ \Top $ for the category of \textit{numerically generated} topological spaces \Cref{nul:terminology}.

\begin{recollection}
	The \defn{Alexandroff topology} on a poset $ P $ is the topology on the underlying set of $ P $ in which a subset $ U \subset P $ is open if and only if $ x \in U $ and $ y \geq x $ implies that $ y \in U $.
\end{recollection}

\begin{nul}\label{nul:posetsnumericallygen}
	Note that every poset in the Alexandroff topology is a numerically generated topological space.


\end{nul}

\begin{definition}\label{def:P-stratified_spaces}
	Let $ P $ be a poset.
	We simply write $ P \in \Top $ for the set $ P $ equipped with the Alexandroff topology.
	The category of \defn{$ P $-stratified topological spaces} is the overcategory $ \TopP $.
	If $ s \colon \fromto{T}{P} $ is a $ P $-stratified topological space, for each $ p \in P $ we write \smash{$ T_p \coloneq s^{-1}(p) $} for the \smash{\defn{$ p $-th stratum}} of $ T $.
\end{definition}

\begin{notation}
	Let $ B $ be a topological space, and $ T,S \in \TopB $.
	We write $ \Map_{/B}(T,S) $ for the topological space of maps $ \fromto{T}{S} $ over $ B $.
	If we need to clarify notation, we write $ \Map_{\TopB}(T,S) $ for this topological space.

	For any topological space $ V $, we write $ T \tensor_B V $ or simply $ T \tensor V $ for the object of $ \TopB $ given by the product $ T \cross V $ with structure morphism induced by the projection $ \fromto{T \cross V}{T} $.
\end{notation}

\begin{nul}\label{nul:Piscolimoversubdiv}
	Let $ P $ be a poset.
	Then since the subdivision $ \sd(P) $ of $ P $ is the category of nondegenerate simplices of $ P $, the poset $ P $ is the colimit $ \colim_{\Sigma \in \sd(P)} \Sigma $ in the category of posets (equivalently, in $ \sSet $).
\end{nul}

\begin{recollection}[{\cite[\HAsec{A.6}]{HA}}]\label{rec:strattopspc}
	Let $ P $ be a poset.
	There is a natural stratification
	\begin{equation*}
		\pi_{P} \colon \fromto{|P|}{P}
	\end{equation*}
	of the geometric realization of (the nerve of) $ P $ by the Alexandroff space $ P $.
	This is defined by appealing to \Cref{nul:Piscolimoversubdiv}, which implies that it suffices to give the standard topological $ n $-simplex $ \real{\Delta^n} $ a $ [n] $-stratification natural in $ [n] $.
	This is given by the map $ \fromto{\real{\Delta^n}}{[n]} $ defined by
	\begin{equation*}
		\goesto{(t_0,\ldots,t_n)}{\max\, \setbar{i \in [n]}{t_i \neq 0}} \period
	\end{equation*}

	If $ X $ is a $ P $-stratified simplicial set, then we can stratify the geometric realization $ \real{X} $ by composing the structure morphism $ \fromto{\real{X}}{\real{P}} $ with $ \pi_P $.
	This defines a left adjoint functor $ \realP{-} \colon \fromto{\sSetP}{\TopP} $ with right adjoint $ \SingP \colon \fromto{\TopP}{\sSetP} $ computed by the pullback of simplicial sets
	\begin{equation*}
		\SingP(T) \colonequals P \crosslimits_{\Sing(P)} \Sing(T) \rlap{\ ,}
	\end{equation*}
	where the morphism $ \fromto{P}{\Sing(P)} $ is adjoint to $ \pi_P $.
\end{recollection}

\begin{nul}\label{nul:stratumofSingP}
	Let $ T $ be a $ P $-stratified topological space.
	Then for each $ p \in P $, the stratum $ \SingP(T)_p $ is isomorphic to the Kan complex $ \Sing(T_p) $.
\end{nul}

\begin{nul}\label{nul:conicstratisqcat}
	Lurie proves \HAa{Theorem}{A.6.4} that if \smash{$ T \in \TopP $} is \textit{conically stratified}\footnote{See \cite[Definitions \HAappthmlink{A.5.3} \& \HAappthmlink{A.5.5}]{HA} for the definition of a conically stratified topological space.}, then the simplicial set \smash{$ \SingP(T) $} is a quasicategory.
\end{nul}

We will use the following observation repeatedly throughout this text.

\begin{remark}\label{rem:Singoftoplink}
	Let $ P $ be a poset.
	Then the adjunction $ \adjto{\realP{-}}{\sSetP}{\TopP}{\SingP} $ is simplicial.
	That is, if $ X $ is a $ P $-stratified simplicial set and $ T $ is a $ P $-stratified topological space, then we have a natural isomorphism of simplicial sets
	\begin{equation*}
		\Sing(\Map_{\TopP}(|X|_P,T)) \isomorphic \Map_{\sSetP}(X,\SingP(T)) \period
	\end{equation*}
\end{remark}


\subsection{Stratified topological spaces as décollages}\label{sec:accessible_loc}

In this subsection we prove \Cref{mainthm:TopPpresentsStrP}.
First we set some notation for the adjunction relating $ P $-stratified topological spaces and simplicial presheaves on the subdivision of $ P $ and recall Douteau's Transfer Theorem (\Cref{thm:Douteaumain}).

\begin{notation}\label{ntn:NPandLP}
	Let $ P $ be a poset.
	We write
	\begin{align*}
		\Nerve_P \colon \sSetP &\to \Fun(\sdop{P},\sSet)
	\intertext{for the functor given by the assignment}
		X &\mapsto [\Sigma \mapsto \Map_{/P}(\Sigma, X)] \period
	\end{align*}
	The functor $ \Nerve_P $ admits a left adjoint $ \LP \colon \fromto{\Fun(\sdop{P},\sSet)}{\sSetP} $ given by the left Kan extension of the Yoneda embedding $ \incto{\sd(P)}{\Fun(\sdop{P},\sSet)} $ along the fully faithful functor $ \incto{\sd(P)}{\sSetP} $ given by $ \goesto{\Sigma}{[\Sigma \subset P]} $.
	Thus $ \LP $ is given by the coend formula\footnote{See \cite[\S1.5]{MR3221774}.}
	\begin{equation*}
		\LP(F) \isomorphic \int^{\Sigma \in \sd(P)} \Sigma \tensor F(\Sigma) \period
	\end{equation*}
\end{notation}

\begin{remark}
	Let $ P $ be a poset.
	Write $ \Pair(P) \subset \sdop{P} \cross \sd(P) $ for the full subposet spanned by those pairs $ (\Sigma,\Sigma') $ where $ \Sigma' \subset \Sigma $.
	The poset $ \Pair(P) $ is an explicit description of the opposite of the twisted arrow category of $ \sdop{P} $.
	Hence by the formula for a coend in terms of a colimit over twisted arrow categories (see \cite[Chapter XI, \S5, Proposition 1]{MR1712872}), the value of the left adjoint $ \LP $ on a functor $ F \colon \fromto{\sdop{P}}{\sSet} $ is given by the colimit
	\begin{equation*}
		\LP(F) \isomorphic \colim_{(\Sigma,\Sigma') \in \Pair(P)} \Sigma' \tensor F(\Sigma) \period
	\end{equation*}
\end{remark}

\begin{notation}\label{ntn:DP}
	Write $ \DP \colon \fromto{\TopP}{\Fun(\sdop{P},\sSet)} $ for the right adjoint functor given by the composite $ \Nerve_P \of \SingP $.
	It follows from \Cref{rem:Singoftoplink} that $ \DP $ is given by the assignment
	\begin{equation*}
		T \mapsto [\Sigma \mapsto \Sing \Map_{\TopP}(\realP{\Sigma},T)] \period
	\end{equation*}
\end{notation}

\begin{notation}
	Let $ P $ be a poset.
	We write $ \sSetPproj $ for category of simplicial presheaves on $ \sd(P) $ given the projective model structure with respect to the Kan model structure on $ \sSet $.
\end{notation}

The following is Douteau's Transfer Theorem.
For the proof, see \cites[Théorèmes 7.2.1, 7.3.7, 7.3.8 \& 7.3.10]{arXiv:1908.01366}[Theorems 2.12 \& 2.15]{MR4224748}[Theorem 3.15]{arXiv:2102.04876}.

\begin{theorem}[(Douteau)]\label{thm:Douteaumain}
	Let $ P $ be a poset.
	The projective model structure on $ \Fun(\sdop{P},\sSet) $ right-transfers to $ \TopP $ along the simplicial adjunction
	\begin{equation}\label{eq:transferadjunction}
		\adjto{\realP{-} \of \LP}{\sSetPproj}{\TopP}{\DP} \period 
	\end{equation}
	Moreover, with respect to these model structures, the adjunction \eqref{eq:transferadjunction} is a simplicial Quillen equivalence of combinatorial simplicial model categories.
\end{theorem}

\begin{warning}
	The proof of \Cref{thm:nerveequiv} that Douteau presented in \cite[Théorèmes 7.2.1, 7.3.7, 7.3.8 \& 7.3.10]{arXiv:1908.01366} and \cite[Theorems 2.12 \& 2.15]{MR4224748} contains a gap. 
	Douteau has since filled in this gap; see \cite[Theorem 3.15]{arXiv:2102.04876}.
\end{warning}

\begin{definition}
	We refer to the model structure on $ \TopP $ of \Cref{thm:Douteaumain} as the \defn{Douteau--Henriques} model structure.
\end{definition}

\begin{nul}\label{nul:explicitdiscofwefibcofib}
	The Douteau--Henriques model structure on $ \TopP $ admits the following explicit description:
	\begin{enumerate}[label=\stlabel{nul:explicitdiscofwefibcofib}, ref=\arabic*]
		\item A morphism $ f \colon \fromto{T}{S} $ in $ \TopP $ is a Douteau--Henriques fibration if and only if for every chain $ \Sigma \subset P $, the induced map of topological spaces
		\begin{equation*}
			\fromto{\MapP(\realP{\Sigma},T)}{\MapP(\realP{\Sigma},S)}
		\end{equation*}
		is a Serre fibration.

		\item A morphism $ f \colon \fromto{T}{S} $ in $ \TopP $ is a Douteau--Henriques weak equivalence if and only if for every chain $ \Sigma \subset P $, the induced map of topological spaces
		\begin{equation*}
			\fromto{\MapP(\realP{\Sigma},T)}{\MapP(\realP{\Sigma},S)}
		\end{equation*}
		is a weak homotopy equivalence.

		\item The sets
		\begin{equation*}
			\left\lbrace\, \incto{\realP{\Sigma \tensor \partial\Delta^n}}{\realP{\Sigma \tensor \Delta^n}} \,\big\vert\, \Sigma \in \sd(P), n \geq 0 \,\right\rbrace
		\end{equation*}
		and 
		\begin{equation*}
			\left\lbrace\, \incto{|\Sigma \tensor \horn{k}{n}|_P}{\realP{\Sigma \tensor \Delta^n}} \,\big\vert\, \Sigma \in \sd(P), n \geq 0 \text{ and } 0 \leq k \leq n \,\right\rbrace
		\end{equation*}
		are generating sets of Douteau--Henriques cofibrations and trivial cofibrations, respectively.
	\end{enumerate}
\end{nul}

\begin{nul}\label{nul:whatTopPpresents}
	The \category $ \Fun(\sdop{P},\Spc) $ of presheaves of \groupoids on $ \sd(P) $ is the underlying \category of the combinatorial simplicial model category $ \sSetPproj $ \HTT{Proposition}{4.2.4.4}.
	Hence the simplicial Quillen equivalence \eqref{eq:transferadjunction} provides an equivalence of \categories between the underlying \category of $ \TopP $ and $ \Fun(\sdop{P},\Spc) $.
\end{nul}

\begin{nul}\label{nul:equivalencesforgoodspaces1}
	\Cref{thm:nerveequiv,rem:Singoftoplink} show that if $ f \colon \fromto{T}{S} $ is a morphism in \smash{$ \TopP $} and both $ \SingP(T) $ and $ \SingP(S) $ are quasicategories, then $ f $ is a Douteau--Henriques equivalence if and only if $ \SingP(f) $ is an equivalence when regarded as a morphism in the \category $ \StrP $ of abstract $ P $-stratified homotopy types. 
\end{nul}

We now arrive at the main result of this section:

\begin{theorem}\label{thm:TopPpresentsStrP}
	Let $ P $ be a poset.
	Then the \category $ \StrP $ is equivalent to an $ \upomega $-accessible localization of the underlying \category of the combinatorial simplicial model category $ \TopP $.
\end{theorem}

\begin{proof}
	Since the underlying \category of $ \TopP $ is equivalent to $ \Fun(\sdop{P},\Spc) $ and $ \StrP $ is an $ \upomega $-accessible localization of $ \Fun(\sdop{P},\Spc) $ (\Cref{thm:nerveequiv,nul:StrPisanaccessibleloc}), we deduce that $ \StrP $ is an $ \upomega $-accessible localization of the underlying \category of $ \TopP $. 
\end{proof}


\begin{nul}\label{nul:leftmodelcat}
	Since the model category $ \sSetPproj $ is \textit{left proper}, there exists a left Bousfield localization of the projective model structure presenting the \category $ \StrP $ \HTT{Proposition}{A.3.7.8}.
	We do not, however, know whether or not the Douteau--Henriques model structure on $ \TopP $ is left proper.
	So while there does exist a left Bousfield localization of $ \TopP $ presenting the \category $ \StrP $, we only know that it exists as a \textit{left model category} \cites[4.13]{MR2771591}{arXiv:2001.03764}[Theorem 7.3]{arXiv:2005.02360}, and it may not exist as a model category.

	In any case, \Cref{thm:TopPpresentsStrP} shows that \category can be obtained from the ordinary category $ \TopP $ of $ P $-stratified topological spaces by inverting a class of weak equivalences (in the \categorical sense).
\end{nul}


\section{The Joyal--Kan model structure}\label{sec:modelonsSetP}

In this section we define a combinatorial simplicial model structure on $ \sSetP $ that presents the \category $ \StrP $.
\Cref{sec:JoyalKandef,subsec:JKfibrant,sec:goodhorns} explore the basic properties of this model structure, and \Cref{sec:simpliciality,subsec:stabfilteredcolim} are dedicated to proving it is simplicial.


\subsection{Definition}\label{sec:JoyalKandef}

In this subsection we define a \textit{Joyal--Kan} model structure on simplicial sets stratified over a poset $ P $ by taking the left Bousfield localization of the Joyal model structure that inverts those simplicial homotopies $ \fromto{X \cross \Delta^1}{Y} $ over $ P $ respecting stratifications.
 
\begin{notation}
	Let $ P $ be a poset.
	Write $ \EP $ for the set of morphisms in $ \sSetP $ consisting of the endpoint inclusions $ \Delta^{\{0\}} , \Delta^{\{1\}} \subset \Delta^1 $ over $ P $ for which the stratification $ f \colon \fromto{\Delta^1}{P} $ is constant.
\end{notation}

\begin{definition}\label{def:Joyal-Kan}
	Let $ P $ be a poset.
	The \defn{Joyal--Kan model structure} on $ \sSetP $ is the $ \sSetJoyal $-enri\-ched left Bousfield localization of the Joyal model structure on $ \sSetP $ with respect to the set $ \EP $. 
\end{definition}

Since the Joyal model structure on $ \sSetP $ is $ \sSetJoyal $-enriched, \cite[Theorems 4.7 \& 4.46]{MR2771591} shows that the Joyal--Kan model structure on $ \sSetP $ exists and satisfies the expected properties:

\begin{proposition}\label{prop:JKexists}
	Let $ P $ be a poset.
	The Joyal--Kan model structure on $ \sSetP $ exists and satisfies the following properties.
	\begin{enumerate}[label=\stlabel{prop:JKexists}, ref=\arabic*]
		\item\label{prop:JKexists.1} The Joyal--Kan model structure on $ \sSetP $ is combinatorial.

		\item\label{prop:JKexists.2} The Joyal--Kan model structure on $ \sSetP $ is $ \sSetJoyal $-enriched.

		\item\label{prop:JKexists.3} The cofibrations in the Joyal--Kan model structure are precisely the monomorphisms of simplicial sets.
		In particular, the Joyal--Kan model structure is left proper.

		\item\label{prop:JKexists.4} The fibrant objects in the Joyal--Kan model structure are precisely the fibrant objects in the Joyal model structure on $ \sSetP $ that are also $ \EP $-local.

		\item\label{prop:JKexists.5} The weak equivalences in the Joyal--Kan model structure are the $ \EP $-local weak equivalences.

		\item\label{prop:JKexists.6} Given a morphism $ f \colon \fromto{X}{Y} $ in $ \sSetP $, if $ X $ and $ Y $ are Joyal--Kan fibrant, then $ f $ is a Joyal--Kan fibration if and only if $ f $ is a Joyal fibration.

		\item\label{prop:JKexists.7} Given a morphism $ f \colon \fromto{X}{Y} $ in $ \sSetP $, if $ X $ and $ Y $ are Joyal--Kan fibrant, then $ f $ is a Joyal--Kan weak equivalence if and only if $ f $ is a Joyal weak equivalence.
	\end{enumerate}
\end{proposition}

\noindent For points \enumref{prop:JKexists}{6} and \enumref{prop:JKexists}{7}, see \cite[Propositions 4.1.7 \& 4.1.8]{Balchin:Model_categories}.

\begin{remark}\label{rem:JKPtrivial}
	When $ P = \ast $ is the terminal poset, the Joyal--Kan model structure on $ \sSet = \sSet_{/*} $ coincides with the Kan model structure.
\end{remark}


\subsection{Fibrant objects in the Joyal–Kan model structure}\label{subsec:JKfibrant}

We now identify the fibrant objects in the Joyal--Kan model structure.

\begin{recollection}\label{2.rec:isofib}
	By \HTT{Corollary}{2.4.6.5} if $ C $ is a quasicategory, then a morphism of simplicial sets $ f \colon \fromto{X}{C} $ is a fibration in the Joyal model structure on $ \sSet $ if and only if the following conditions are satisfied:
	\begin{enumerate}[label=\stlabel{2.rec:isofib}, ref=\arabic*]
		\item\label{2.rec:isofib.1} The morphism $ f $ is an inner fibration.

		\item\label{2.rec:isofib.2} For every equivalence $ e \colon \equivto{c}{c'} $ in $ C $ and object $ \ctilde \in X $ such that $ f(\ctilde) = c $, there exists an equivalence $ \etilde \colon \equivto{\ctilde}{\ctilde'} $ in $ X $ such that $ f(\etilde) = e $.
	\end{enumerate}
	A morphism of simplicial sets satisfying \enumref{2.rec:isofib}{1} and \enumref{2.rec:isofib}{2} is called an \defn{isofibration}.
	(See also \cite[\S 2]{MR3792514}.)
\end{recollection}

\noindent We make use of the following immediate consequence of the definitions.

\begin{lemma}\label{lem:isofibtoposet}
	Let $ C $ be a quasicategory whose equivalences are precisely the degenerate edges (e.g., a poset).
	Then a morphism of simplicial sets $ f \colon \fromto{X}{C} $ is an isofibration if and only if $ f $ is an inner fibration.
\end{lemma}

\begin{proposition}\label{prop:identifyfibrant}
	Let $ P $ be a poset.
	An object $ X $ of $ \sSetP $ is fibrant in the Joyal--Kan model structure if and only if the structure morphism $ \fromto{X}{P} $ is an inner fibration and for every $ p \in P $ the stratum $ X_p $ is a Kan complex.
\end{proposition}

\begin{proof}
	Since the Joyal--Kan model structure on $ \sSetP $ is the left Bousfield localization of the Joyal model structure on $ \sSetP $ with respect to $ \EP $, the fibrant objects in the Joyal--Kan model structure on $ \sSetP $ are the fibrant objects in the Joyal model structure on $ \sSetP $ that are also $ \EP $-local.
	An object $ X \in \sSetP $ is fibrant in the Joyal model structure if and only if the structure morphism $ \fromto{X}{P} $ is an isofibration, or, equivalently the structure morphism $ \fromto{X}{P} $ is an inner fibration (\Cref{lem:isofibtoposet}).

	Now we analyze the $ \EP $-locality condition.
	A Joyal-fibrant object $ X \in \sSetP $ is $ \EP $-local if and only if for every $ 1 $-simplex \smash{$ \sigma \colon \fromto{\Delta^1}{P} $} such that $ \sigma(0) = \sigma(1) $, evaluation morphisms
	\begin{equation*}
		\ev_i \colon \fromto{\MapP(\Delta^1,X)}{\MapP(\Delta^{\{i\}},X)} 
	\end{equation*}
	for $ i = 0,1 $ are isomorphisms in the homotopy category of $ \sSetJoyal $.
	Let $ p \in P $ be such that $ \sigma(0) = \sigma(1) = p $.
	Then
	\begin{equation*}
		\MapP(\Delta^1,X) \isomorphic \Map(\Delta^1,X_p)
	\end{equation*}
	and
	\begin{equation*}
		\MapP(\Delta^{\{i\}},X) \isomorphic \Map(\Delta^{\{i\}},X_p) \comma
	\end{equation*}
	for $ i=0,1 $.
	Under these identifications, the evaluation morphisms
	\begin{equation*}
		\ev_i \colon \fromto{\MapP(\Delta^1,X)}{\MapP(\Delta^{\{i\}},X)} 
	\end{equation*}
	are identified with the evaluation morphisms
	\begin{equation*}
		\ev_i \colon \fromto{\Map(\Delta^1,X_p)}{\Map(\Delta^{\{i\}},X_p) \isomorphic X_p }  \comma
	\end{equation*}
	for $ i = 0,1 $.
	Since the strata of $ X $ are quasicategories, $ X $ is $ \EP $-local if and only if for every $ p \in P $, the evaluation morphisms
	\begin{equation*}
		\ev_i \colon \fromto{\Map(\Delta^1,X_p)}{\Map(\Delta^{\{i\}},X_p) \isomorphic X_p}  \comma
	\end{equation*}
	for $ i = 0,1 $, are Joyal equivalences. 
	To conclude, recall that for a quasicategory $ C $, the evaluation morphisms $ \ev_0,\ev_1 \colon \fromto{\Map(\Delta^1,C)}{C} $ are Joyal equivalences if and only if $ C $ is a Kan complex.
\end{proof}

Combining \Cref{prop:identifyfibrant} with \HTT{Proposition}{2.3.1.5} we deduce:

\begin{proposition}\label{prop:equivfibrant}
	Let $ P $ be a poset, $ X $ a simplicial set, and $ f \colon \fromto{X}{P} $ a morphism of simplicial sets.
	The following are equivalent:
	\begin{enumerate}[label=\stlabel{prop:equivfibrant}, ref=\arabic*]
		\item\label{prop:equivfibrant.1} The object $ f \colon \fromto{X}{P} $ of $ \sSetP $ is fibrant in the Joyal--Kan model structure.

		\item\label{prop:equivfibrant.2} The morphism $ f \colon \fromto{X}{P} $ is an inner fibration with all fibers Kan complexes.

		\item\label{prop:equivfibrant.3} The simplicial set $ X $ is a quasicategory and all of the fibers of $ f $ are Kan complexes.

		\item\label{prop:equivfibrant.4} The simplicial set $ X $ is a quasicategory and $ f $ is a conservative functor between quasicategories.
	\end{enumerate}
\end{proposition}




\begin{corollary}\label{cor:morfibiscons}
	Let $ P $ be a poset.
	A morphism in $ \sSetP $ between fibrant objects in the Joyal--Kan model structure is a conservative functor.
\end{corollary}

\begin{proof}
	Note that if a composite functor $ gf $ is conservative and $ g $ is conservative, then $ f $ is conservative.
\end{proof}


\subsection{Stratified horn inclusions}\label{sec:goodhorns}

In this subsection we characterize the horn inclusions in $ \sSetP $ that are Joyal--Kan equivalences.
We will use these horn inclusions in our proof that the Joyal--Kan model structure is simplicial (see \cref{sec:simpliciality}).

\begin{proposition}\label{prop:stratifiedhornsall}
	Let $ P $ be a poset.
	A horn inclusion $i\colon\incto{\Lambda^n_k}{\Delta^n}$ over $ P $ stratified by a morphism $f\colon\fromto{\Delta^n}{P}$ is a Joyal--Kan equivalence in $ \sSetP $ if and only if one of the following conditions holds:
	\begin{enumerate}[label=\stlabel{prop:stratifiedhornsall}, ref=\arabic*]
		\item\label{prop:stratifiedhornsall.1} $ 0<k<n $.

		\item\label{prop:stratifiedhornsall.2} $ k=0 $ and $ f(0) = f(1) $.

		\item\label{prop:stratifiedhornsall.3} $ k=n $ and $ f(n-1) = f(n) $.
	\end{enumerate}
\end{proposition}

\begin{proof}
	First we show that the class of horn inclusions \enumref{prop:stratifiedhornsall}{1}--\enumref{prop:stratifiedhornsall}{3} are Joyal--Kan equivalences.
	It is clear that inner horn inclusions $ \incto{\horn{k}{n}}{\Delta^n} $ are weak equivalences in the Joyal--Kan model structure on $ \sSetP $ as they are already weak equivalences in the Joyal model structure.
	If $ n = 1 $, then the endpoint inclusions $ \incto{\horn{0}{1},\horn{1}{1}}{\Delta^1} $ where $ f(0) = f(1) $ are Joyal--Kan equivalences by the definition of the Joyal--Kan model structure.

	Now we treat the case of higher outer horns.
	We treat the case of left horns $ \incto{\horn{0}{n}}{\Delta^n} $ where the stratification $ f \colon \fromto{\Delta^n}{P} $ has the property that $ f(0) = f(1) $ (i.e., the class specified by \enumref{prop:stratifiedhornsall}{2}); the case of right horns is dual.
	We prove the claim by induction on $ n $.

	For the base case where $ n = 2 $, write $ D_0^2 $ for the (nerve of the) preorder given by $ 0 \leq 1 \leq 2 $ along with $ 0 \geq 1 $, and stratify $ D_0^2 $ by the unique extension of $ f $ to $ D_0^2 $.
	All stratifications will be induced by $ f $ via inclusions into $ D_0^2 $.
	We prove the claim by showing that the inclusions $ \incto{\horn{0}{2},\Delta^2}{D_0^2} $ are Joyal--Kan equivalences and conclude by the 2-of-3 property.
	Write $ E $ for the walking isomorphism category $ 0 \isomorphic 1 $ and consider the cube
	\begin{equation}\label{2.cube:pushout}
		\begin{tikzcd}[column sep={10ex,between origins}, row sep={8ex,between origins}]
			& \Delta^{\{0\}} \arrow[dl, equals] \arrow[rr, hooked] \arrow[dd, hooked, "\wr"{xshift=-0.1em, near start}] \arrow[ddrr, phantom, very near end, "\ulcorner", yshift=-0.2em, xshift=0.6em] & & \Delta^{\{0<2\}} \arrow[dl, equals] \arrow[dd, hooked, "\wr"{xshift=-0.1em}] \\
			\Delta^{\{0\}} \arrow[dd, hooked, "\wr"{xshift=-0.1em}] \arrow[rr, crossing over, hooked] \arrow[ddrr, phantom, very near end, "\ulcorner", yshift=-0.2em, xshift=0.6em] & & \Delta^{\{0<2\}} \\
			& \Delta^{\{0<1\}} \arrow[dl, hooked'] \arrow[rr, hooked] & & \horn{0}{2} \arrow[dl, hooked'] \\
			E \arrow[rr, hooked] & & L_0^2 \arrow[from=uu, crossing over, hooked, "\wr"{xshift=-0.1em, near end}] & \phantom{\horn{0}{2}} \comma
		\end{tikzcd}
	\end{equation}
	where the front face is a pushout defining the simplicial set $ L_0^2 $ and the back face is a pushout square.
	Since $ f(0) = f(1) $, the inclusion $ \incto{\Delta^{\{0\}}}{\Delta^{\{0<1\}}} $ is a trivial Joyal--Kan cofibration; the fact that the back face of \eqref{2.cube:pushout} is a pushout then shows that the inclusion $ \incto{\Delta^{\{0<2\}}}{\horn{0}{2}} $ is a trivial Joyal--Kan cofibration.
    Since the inclusion $ \incto{\Delta^{\{0\}}}{E} $ is a trivial Joyal cofibration and the front face of \eqref{2.cube:pushout} is a pushout, the inclusion $ \incto{\Delta^{\{0<2\}}}{L_0^2} $ is a trivial Joyal cofibration.
    By the 2-of-3 property, the induced map on pushouts $ \incto{\horn{0}{2}}{L_0^2} $ is a trivial Joyal--Kan cofibration.
    Similarly, the inclusion
    \begin{equation*}
    	\incto{\horn{1}{2}}{L_1^2 \colonequals E \union^{\Delta^{\{1\}}} \Delta^{\{1<2\}}} 
    \end{equation*}
    is a trivial Joyal--Kan cofibration.
    The inclusions $ \incto{L_0^2, L_{1}^2}{D_0^2} $ are trivial Joyal cofibrations, so in particular the composite inclusion
    \begin{equation*}
    	\incto{\horn{0}{2}}{\incto{L_0^2}{D_0^2}}
    \end{equation*}
    is a trivial Joyal--Kan cofibration.
    Finally, to see that the inclusion $ \incto{\Delta^2}{D_0^2} $ is a Joyal--Kan equivalence note that we have a commutative square
    \begin{equation*}
      \begin{tikzcd}
       \horn{1}{2} \arrow[d, hooked, "\wr"{xshift=-0.1em}] \arrow[r, hooked, "\sim"{yshift=-0.2em}] & \Delta^{2} \arrow[d, hooked] \\ 
       L_1^2 \arrow[r, hooked, "\sim"{yshift=-0.2em}] & D_0^2 \comma
      \end{tikzcd}
    \end{equation*}
    where the horizontal morphisms are trivial Joyal cofibrations and the inclusion $ \incto{\horn{1}{2}}{L_1^2} $ is a Joyal--Kan equivalence.
    This concludes the base case.

	Now we prove the induction step with $ n \geq 3 $ and $ \incto{\horn{0}{n}}{\Delta^n} $ an outer horn inclusion over $ P $ where the stratification $ f \colon \fromto{\Delta^n}{P} $ has the property that $ f(0) = f(1) $.
	Write
	\begin{equation*}
		\horn{\{0<2\}}{n} \colonequals \bigcup_{j \in [n] \sminus \{0<2\}} \Delta^{\{0<\cdots<j-1,j+1<\cdots<n\}} \subset \Delta^n 
	\end{equation*}
	and note that by \cite[Lemma 12.13]{MR3558219} the inclusion $ \incto{\horn{\{0<2\}}{n}}{\Delta^n} $ is inner anodyne.
	Since we have a factorization of the inclusion $ \incto{\horn{\{0<2\}}{n}}{\Delta^n} $ as a composite
	\begin{equation*}
		\incto{\incto{\horn{\{0<2\}}{n}}{\horn{0}{n}}}{\Delta^n} \comma 
	\end{equation*}
	the claim is equivalent to showing that the inclusion $ \incto{\horn{\{0<2\}}{n}}{\horn{0}{n}} $ is a Joyal--Kan equivalence in $ \sSetP $.
	To see this, note that we have a pushout square in $ \sSetP $ 
	\begin{equation}\label{2.sq:noliftspaces}
      \begin{tikzcd}
       \horn{0}{\{0<1<3<\cdots<n\}}  \arrow[dr, phantom, very near end, "\ulcorner", xshift=0.4em, yshift=-0.1em] \arrow[d, hooked] \arrow[r, hooked, "\sim"{yshift=-0.2em}] & \Delta^{\{0<1<3<\cdots<n\}} \arrow[d, hooked] \\ 
       \horn{\{0<2\}}{n} \arrow[r, hooked, "\sim"{yshift=-0.2em}] & \horn{0}{n} \comma
      \end{tikzcd}
    \end{equation}
    where the inclusion $ \incto{\horn{0}{\{0<1<3<\cdots<n\}}}{\Delta^{\{0<1<3<\cdots<n\}}} $ is a trivial Joyal--Kan cofibration by the induction hypothesis.

	Now we prove the horn inclusions given by the classes \enumref{prop:stratifiedhornsall}{1}--\enumref{prop:stratifiedhornsall}{3} are the only horn inclusions over $ P $ that are trivial Joyal--Kan cofibrations.
	Equivalently, if $ i \colon \incto{\horn{k}{n}}{\Delta^n} $ is an outer horn and either $ k = 0 $ and $ f(0) \neq f(1) $, or $ k = n $ and $ f(n-1) \neq f(n) $, then $ i $ is not a Joyal--Kan equivalence.
	We treat the case that $ k = 0 $; the case that $ k = n $ is dual.
	The cases where $ n = 1 $ and $ n = 2 $ require slightly different (but easier) arguments than when $ n \geq 3 $, so we treat those first.

	When $ n = 1 $, we need to show that the endpoint inclusion $ \incto{\Delta^{\{0\}}}{\Delta^1} $ is not a Joyal--Kan equivalence, where the stratification $ f \colon \fromto{\Delta^1}{P} $ is a monomorphism.
	In this case, by \Cref{prop:equivfibrant} both $ \Delta^{\{0\}} $ and $ \Delta^1 $ are fibrant in the Joyal--Kan model structure, so by \enumref{prop:JKexists}{7} we just need to check that the inclusion $ i \colon \incto{\Delta^{\{0\}}}{\Delta^1} $ is not a Joyal equivalence, which is clear.

	For $ n = 2 $, note that the simplicial set $ \horn{0}{2} $ is a $ 1 $-category.
	Since $ f(0) \neq f(1) $, we have $ f(0) \neq f(2) $, so the functor $ f \colon \fromto{\horn{0}{2}}{P} $ is conservative; applying \Cref{prop:equivfibrant} shows that $ \horn{0}{2} $ is fibrant in the Joyal--Kan model structure.
	To see that the inclusion $ \incto{\horn{0}{2}}{\Delta^2} $ is not a trivial Joyal--Kan cofibration, note that the lifting problem
	\begin{equation*}
      \begin{tikzcd}
      	\horn{0}{2} \arrow[d, hooked] \arrow[r, equals] & \horn{0}{2} \arrow[d, ->>, "f"] \\ 
      	\Delta^2 \arrow[r, "f"'] \arrow[ur, dotted] & P 
      \end{tikzcd}
    \end{equation*}
    does not admit a solution because the inclusion of simplicial sets $ \incto{\horn{0}{2}}{\Delta^2} $ does not admit a retraction.

    To prove the claim for $ n \geq 3 $, one can easily construct a $ 1 $-category $ C_{0,f}^n $ along with a natural inclusion $ \phi_f \colon \incto{\horn{0}{n}}{C_{0,f}^n} $ that does not extend to $ \Delta^n $ as follows: adjoin a new morphism $ a \colon \fromto{1}{n} $ to $ \Delta^n $ so that $ a $ and the unique morphism $ \fromto{1}{n} $ are equalized by the unique morphism $ \fromto{0}{1} $, then formally adjoin inverses to all morphisms $ \fromto{i}{j} $ such that $ f(i) = f(j) $.
    The inclusion $ \phi_f \colon \incto{\horn{0}{n}}{C_{0,f}^n} $ is not the standard one, but one with the property that the edge $ \Delta^{\{1<n\}} $ is sent to the morphism $ a $.
    Thus $ \phi_f $ does not extend to $ \Delta^n $.
    The morphism $ \restrict{f}{\horn{0}{n}} $ extends to a stratification $ \fbar \colon \fromto{C_{0,f}^n}{P} $ that makes $ C_{0,f}^n $ a fibrant object in the Joyal--Kan model structure, and the inclusion $ \incto{\horn{0}{n}}{\Delta^n} $ is not a trivial Joyal--Kan cofibration in $ \sSetP $ since the lifting problem  
    \begin{equation*}
      \begin{tikzcd}
       \horn{0}{n} \arrow[d, hooked] \arrow[r, hooked, "\phi_f"] & C_{0,f}^n \arrow[d, ->>, "\fbar"] \\ 
       \Delta^n \arrow[r, "f"'] \arrow[ur, dotted] & P 
      \end{tikzcd}
    \end{equation*}
    does not admit a solution.
\end{proof}

\begin{notation}
	Let $ P $ be a poset.
	Write $ \JP \subset \Mor(\sSetP) $ for the set of all horn inclusions $ i \colon \incto{\Lambda^n_k}{\Delta^n}$ over $ P $ that are Joyal--Kan equivalences.
\end{notation}

We can use the set $ \JP $ to identify fibrations between fibrant objects of the Joyal--Kan model structure on $ \sSetP $.
First we record a convenient fact.

\begin{lemma}\label{lem:consisofib}
	Let $ f \colon \fromto{X}{Y} $ be a conservative functor between quasicategories.
	The following are equivalent: 
	\begin{enumerate}[label=\stlabel{lem:consisofib}, ref=\arabic*]
		\item\label{lem:consisofib.1} For every equivalence $ e \colon \equivto{y}{y'} $ in $ Y $ and object $ \ytilde \in X $ such that $ f(\ytilde) = y $, there exists an equivalence $ \etilde \colon \equivto{\ytilde}{\ytilde'} $ in $ X $ such that $ f(\etilde) = e $.

		\item\label{lem:consisofib.2} For every equivalence $ e \colon \equivto{y}{y'} $ in $ Y $ and object $ \ytilde \in X $ such that $ f(\ytilde) = y $, there exists a morphism $ \etilde \colon \fromto{\ytilde}{\ytilde'} $ in $ X $ such that $ f(\etilde) = e $.
	\end{enumerate}
\end{lemma}


\begin{proposition}\label{prop:fibfibrant}
	Let $ P $ be a poset and $ f \colon \fromto{X}{Y} $ a morphism in $ \sSetP $ between fibrant objects in the Joyal--Kan model structure.
	Then the following are equivalent:
	\begin{enumerate}[label=\stlabel{prop:fibfibrant}, ref=\arabic*]
		\item\label{prop:fibfibrant.1} The morphism $ f $ is a Joyal--Kan fibration.

		\item\label{prop:fibfibrant.2} The morphism $ f $ is a Joyal fibration, equivalently, an isofibration.

		\item\label{prop:fibfibrant.3} The morphism $ f $ satisfies the right lifting property with respect to $ \JP $.

		\item\label{prop:fibfibrant.4} The morphism $ f $ is an inner fibration and the restriction of $ f $ to each stratum is a Kan fibration.

		\item\label{prop:fibfibrant.5} The morphism $ f $ is an inner fibration and satisfies the right lifting property with respect to $ \EP $.
	\end{enumerate}
\end{proposition}

\begin{proof}
	The equivalence \enumref{prop:fibfibrant}{1}$ \Leftrightarrow $\enumref{prop:fibfibrant}{2} is the content of \enumref{prop:JKexists}{6}.

	Now we show that \enumref{prop:fibfibrant}{2}$ \Rightarrow $\enumref{prop:fibfibrant}{3}.
	Assume that $ f $ is an isofibration.
    Since $ f $ is an isofibration, $ f $ is an inner fibration, hence lifts against \textit{inner} horns in $ \JP $.
	Now consider the lifting problem
	\begin{equation}\label{2.diag:edgelifting}
    	\begin{tikzcd}
	       \Delta^{\{i\}} \arrow[d, hooked] \arrow[r, "h"] & X \arrow[d, "f"] \\ 
	       \Delta^1 \arrow[r, "h'"'] \arrow[ur, dotted] & Y  
	    \end{tikzcd}
    \end{equation}
    where the inclusion $ \incto{\Delta^{\{i\}}}{\Delta^1} $ is in $ \JP $.
    Since $ Y $ is fibrant in the Joyal--Kan model structure, the edge $ h'(\Delta^1) $ is an equivalence in $ Y $. 
    \Cref{lem:consisofib} (and its dual) now shows that the lifting problem \eqref{2.diag:edgelifting} admits a solution.
    Finally, if $ n \geq 2 $ and $ k = 0 $ or $ k = n $, then given a lifting problem 
    \begin{equation*}
    	\begin{tikzcd}
	       \horn{k}{n} \arrow[d, hooked] \arrow[r, "h"] & X \arrow[d, "f"] \\
	       \Delta^n \arrow[r, "h'"'] \arrow[ur, dotted] & Y  
	    \end{tikzcd}
    \end{equation*}
    where the horn inclusion $ \incto{\horn{k}{n}}{\Delta^n} $ is in $ \JP $, since $ X $ and $ Y $ are fibrant in the Joyal--Kan model structure:
   	\begin{enumerate}
   		\item If $ k = 0 $, then $ h(\Delta^{\{0<1\}}) $ and $ h'(\Delta^{\{0<1\}}) $ are equivalences.

   		\item If $ k = n $, then $ h(\Delta^{\{n-1<n\}}) $ and $ h'(\Delta^{\{n-1<n\}}) $ are equivalences.
   	\end{enumerate}
   	In either case, the desired lift exists because $ f $ is an inner fibration and the outer horn is ``special'' \cites[Theorem 2.2]{MR1935979}[p. 236]{MR3221774}.

	The fact that \enumref{prop:fibfibrant}{3} implies \enumref{prop:fibfibrant}{4} is clear from the identification of $ \JP $ (\Cref{prop:stratifiedhornsall}).

	The fact that \enumref{prop:fibfibrant}{4} implies \enumref{prop:fibfibrant}{5} is clear from the definition of $ \EP $ and the fact that the restriction of $ f $ to each stratum is a Kan fibration.

	Now we show that \enumref{prop:fibfibrant}{5} implies \enumref{prop:fibfibrant}{2}.
	Assume that $ f $ is an inner fibration and satisfies the right lifting property with respect to $ \EP $.
	Since $ f $ is conservative (\Cref{cor:morfibiscons}) and the equivalences in $ Y $ lie in individual strata, \Cref{lem:consisofib} combined with the fact that $ f $ satisfies the right lifting property with respect to $ \EP $ show that $ f $ is an isofibration.
\end{proof}

\begin{corollary}\label{lem:JPdetectsfibrant}
	Let $ P $ be a poset and $ X $ an object of $ \sSetP $.
	Then $ X $ is fibrant in the Joyal--Kan model structure if and only if the stratification $ \fromto{X}{P} $ satisfies the right lifting property with respect to $ \JP $.
\end{corollary}


\subsection{Simpliciality of the Joyal--Kan model structure}\label{sec:simpliciality}

Unlike the Kan model structure on $ \sSet $, the Joyal model structure is \textit{not} simplicial.
As a result, it does not follow formally from the definition that the Joyal--Kan model structure on $ \sSetP $ is simplicial.
In this subsection we recall three criteria that guarantee that a model structure is simplicial, and verify the first two of them.
We verify the third in \cref{subsec:stabfilteredcolim}.

\begin{remark}\label{2.rmk:pushout-prod}
	Let $ P $ be a poset, $ i \colon \fromto{X}{Y} $ a morphism in $ \sSetP $, and $ j \colon \fromto{A}{B} $ a morphism of simplicial sets.
	Then on underlying simplicial sets, the pushout-tensor
	\begin{equation*}
		i \tensorhat j \colon \fromto{(X \tensor B) \coproductlimits^{X \tensor A} (Y \tensor A)}{Y \tensor B}
	\end{equation*}
	is simply the pushout-product
	\begin{equation*}
		i \crosshat j \colon \fromto{(X \cross B) \coproductlimits^{X \cross A} (Y \cross A)}{Y \cross B}
	\end{equation*}
	in $ \sSet $.
	Since the pushout-product of monomorphisms in $ \sSet $ is a monomorphism and the forgetful functor $ \fromto{\sSetP}{\sSet} $ detects monomorphisms, if $ i $ and $ j $ are monomorphisms, then $ i \tensorhat j $ is a monomorphism.
\end{remark}

\begin{nul}\label{nul:wantsimplcial}
	By appealing to \HTT{Proposition}{A.3.1.7}, we can prove that the Joyal--Kan model structure is simplicial by proving the following three claims:
	\begin{enumerate}[label=\stlabel{nul:wantsimplcial}, ref=\arabic*]
		\item\label{nul:wantsimplcial.1} Given a monomorphism of simplicial sets $ j \colon \into{A}{B} $ and a Joyal--Kan cofibration $ i \colon \into{X}{Y} $ in $ \sSetP $, the pushout-tensor
		\begin{equation*}
			i \tensorhat j \colon \fromto{(X \tensor B) \coproductlimits^{X \tensor A} (Y \tensor A)}{Y \tensor B}
		\end{equation*}
		is a Joyal--Kan cofibration.

		\item\label{nul:wantsimplcial.2} For every $ n \geq 0 $ and every object $ X \in \sSetP $, the natural map 
		\begin{equation*}
			\fromto{X \tensor \Delta^n}{X \tensor \Delta^0 \isomorphic X}
		\end{equation*}
		is a Joyal--Kan equivalence.

		\item\label{nul:wantsimplcial.3} The collection of weak equivalences in the Joyal--Kan model structure on $ \sSetP $ is stable under filtered colimits.
	\end{enumerate}
	Note that \enumref{nul:wantsimplcial}{1} follows from \Cref{2.rmk:pushout-prod} and the fact that cofibrations in the Joyal--Kan model structure are monomorphisms of simplicial sets (\Cref{prop:JKexists}).
\end{nul}

We first concern ourselves with \enumref{nul:wantsimplcial}{2}.
Since the natural map
\begin{equation*}
	\fromto{X \tensor \Delta^n}{X \tensor \Delta^0 \isomorphic X}
\end{equation*}
admits a section $ \incto{X \isomorphic X \tensor \Delta^{\{ 0 \}}}{X \tensor \Delta^n} $, it suffices to show that this section is a trivial Joyal--Kan cofibration.
In fact, we prove a more precise claim.

\begin{notation}\label{ntn:IH_and_LH}
	Let $ P $ be a poset.
	\begin{enumerate}[label=\stlabel{ntn:IH_and_LH}, ref=\arabic*]
		\item Write $ \IH_P \subset \JP $ for the \textit{inner} horn inclusions in $ \JP $.

		\item Write $ \LH_P \subset \JP $ for those horn inclusions $ \incto{\horn{k}{n}}{\Delta^n} $ in $ \JP $ where $ n \geq 1 $ and $ 0 \leq k < n $, i.e., the \textit{left} horn inclusions in $ \JP $.

	\end{enumerate}
	Note that \Cref{prop:stratifiedhornsall} gives complete characterizations of $ \IH_P $ and $ \LH_P $.
\end{notation}

\begin{proposition}\label{prop:incweakequiv}
	Let $ P $ be a poset and $ n \geq 0 $ an integer.
	For any object $ X \in \sSetP $, the inclusion
	\begin{equation*}
		\incto{X \isomorphic X \tensor \Delta^{\{ 0 \}}}{X \tensor \Delta^n}
	\end{equation*}
	is in the weakly saturated class generated by $ \LH_P $, in particular, a trivial Joyal--Kan cofibration in $ \sSetP $.
\end{proposition}

The next proposition (and its proof) is a stratified variant of \cite[Propositions \HTTthmlink{2.1.2.6} \& \HTTthmlink{3.1.1.5}]{HTT}.
We use it to prove \Cref{prop:incweakequiv}.

\begin{proposition}\label{prop:weaklysaturated1}
	Let $ P $ be a poset.
	Consider the following classes of morphisms in $ \sSetP $:
	\begin{enumerate}[label=\stlabel{prop:weaklysaturated1}, ref=\arabic*]
		\item\label{prop:weaklysaturated1.1} All inclusions 
		\begin{equation*}
			\incto{(\partial \Delta^m \tensor \Delta^1) \coproductlimits^{\partial \Delta^m \tensor \Delta^{\{0\}}} (\Delta^m \tensor \Delta^{\{0\}})}{\Delta^m \tensor \Delta^1} \comma
		\end{equation*}
		where $ m \geq 0 $ and $ \Delta^m \in \sSetP $ is any $ m $-simplex over $ P $.

		\item\label{prop:weaklysaturated1.2} All inclusions 
		\begin{equation*}
			\incto{(A \tensor \Delta^1) \coproductlimits^{A \tensor \Delta^{\{0\}}} (B \tensor \Delta^{\{0\}})}{B \tensor \Delta^1} \comma
		\end{equation*}
		where $ \incto{A}{B} $ is any monomorphism in $ \sSetP $.
	\end{enumerate}

	The classes \enumref{prop:weaklysaturated1}{1} and \enumref{prop:weaklysaturated1}{2} generate the same weakly saturated class of morphisms in $ \sSetP $.
	Moreover, this weakly saturated class of morphisms generated by \enumref{prop:weaklysaturated1}{1} or \enumref{prop:weaklysaturated1}{2} is contained in the weakly saturated class of morphisms generated by $ \LH_P $.
\end{proposition}

\begin{proof}
	Since the inclusions $ \incto{\partial \Delta^m}{\Delta^m} $ in $ \sSetP $ generate the monomorphisms in $ \sSetP $, to see that each of the morphisms specified in \enumref{prop:weaklysaturated1}{2} is contained in the weakly saturated class generated by \enumref{prop:weaklysaturated1}{1}, it suffices to work simplex-by-simplex with the inclusion $ \incto{A}{B} $.
	The converse is clear since the class specified by \enumref{prop:weaklysaturated1}{1} is contained in the class specified by \enumref{prop:weaklysaturated1}{2}.

	To complete the proof, we show that for each $ P $-stratified $ m $-simplex $ \Delta^m \in \sSetP $, the inclusion 
	\begin{equation}\label{eq:boundarycylinc}
		\incto{(\partial \Delta^m \tensor \Delta^1) \coproductlimits^{\partial \Delta^m \tensor \Delta^{\{0\}}} (\Delta^m \tensor \Delta^{\{0\}})}{\Delta^m \tensor \Delta^1}
	\end{equation}
	belongs to the weakly saturated class generated by $ \LH_P $.
	The proof of this is \textit{verbatim} the same as the proof of \HTT{Proposition}{2.1.2.6}, which writes the inclusion \eqref{eq:boundarycylinc} as a composite of pushouts of horn inclusions, all of which are in $ \LH_P $.
\end{proof}

\begin{corollary}\label{cor:Delta1inc}
	Let $ P $ be a poset.
	For any $ P $-stratified simplicial set $ X \in \sSetP $, the inclusion $ \incto{X \tensor \Delta^{\{0\}}}{X \tensor \Delta^1} $ is in the weakly saturated class generated by $ \LH_P $.
\end{corollary}

\begin{proof}
	In \enumref{prop:weaklysaturated1}{2}, set $ A = \varnothing $ and $ B = X $.
\end{proof}

\begin{notation}
	Let $ n \geq 0 $ be an integer.
	Write $ \spine{n} \subset \Delta^n $ for the \defn{spine} of $ \Delta^n $, defined by
	\begin{equation*}
		\spine{n} \colonequals \Delta^{\{0<1\}} \union^{\Delta^{\{1\}}} \cdots \union^{\Delta^{\{n-1\}}} \Delta^{\{n-1<n\}} \period
	\end{equation*}
\end{notation}

Now we use \Cref{cor:Delta1inc} and the fact that the spine inclusion $ \incto{\spine{n}}{\Delta^n} $ is inner anodyne to address \Cref{prop:incweakequiv}.

\begin{lemma}\label{lem:spineinc}
	Let $ P $ be a poset and $ n \geq 0 $ an integer.
	For any $ P $-stratified simplicial set $ X \in \sSetP $, the inclusion $ \incto{X \tensor \Delta^{\{0\}}}{X \tensor \spine{n}} $ is in the weakly saturated class generated by $ \LH_P $.
\end{lemma}

\begin{proof}
	Noting that $ \spine{1} = \Delta^1 $, factor the inclusion $ \incto{X \tensor \Delta^{\{0\}}}{X \tensor \spine{n}} $ as a composite
	\begin{equation*}
      \begin{tikzcd}[sep=1.5em]
       X \tensor \Delta^{\{0\}} \arrow[r, hooked] & X \tensor \Delta^1 \arrow[r, hooked] & X \tensor \spine{2} \arrow[r, hooked] & \cdots  \arrow[r, hooked] & X \tensor \spine{n} \period
      \end{tikzcd}
    \end{equation*}
    The inclusion $ \incto{X \tensor \Delta^{\{0\}}}{X \tensor \Delta^1} $ is in the weakly saturated class generated by $ \LH_P $ (\Cref{cor:Delta1inc}), so it suffices to show that for $ 1 \leq k \leq n-1 $, the inclusion $ \incto{X \tensor \spine{k}}{X \tensor \spine{k+1}} $ is in the weakly saturated class generated by $ \LH_P $.
    To see this, note that the inclusion $ \incto{X \tensor \spine{k}}{X \tensor \spine{k+1}} $ is given by the pushout
    \begin{equation*}
      \begin{tikzcd}
      	X \tensor \Delta^{\{k \}}  \arrow[dr, phantom, very near end, "\ulcorner", xshift=0.75em, yshift=-0.25em] \arrow[d, hooked] \arrow[r, hooked] & X \tensor \Delta^{\{k<k+1\}} \arrow[d, hooked] \\ 
       X \tensor \spine{k} \arrow[r, hooked] & X \tensor \spine{k+1} \comma
      \end{tikzcd}
    \end{equation*}
    and by \Cref{cor:Delta1inc} the inclusion $ \incto{X \tensor \Delta^{\{k \}}}{X \tensor \Delta^{\{k<k+1 \}}} $ is in the weakly saturated class generated by $ \LH_P $.
\end{proof}

\begin{proof}[Proof of \Cref{prop:incweakequiv}]
	The inclusion $ \incto{X \tensor \Delta^{\{0\}}}{X \tensor \Delta^n} $ factors as a composite
	\begin{equation*}
      \begin{tikzcd}[sep=1.5em]
       X \tensor \Delta^{\{0\}} \arrow[r, hooked] & X \tensor \spine{n} \arrow[r, hooked] & X \tensor \Delta^n \period
      \end{tikzcd}
    \end{equation*}
    To conclude, first note that by \Cref{lem:spineinc} the inclusion $ \incto{X \tensor \Delta^{\{0\}}}{X \tensor \spine{n}} $ is in the weakly saturated class generated by $ \LH_P $.
    Second, since the inclusion $ \incto{\spine{n}}{\Delta^n} $ is inner anodyne, the inclusion 
    \begin{equation*}
      \begin{tikzcd}[sep=1.5em]
       X \cross \spine{n} = X \tensor \spine{n} \arrow[r, hooked] & X \tensor \Delta^n = X \cross \Delta^n 
      \end{tikzcd}
    \end{equation*}
    is inner anodyne \HTT{Corollary}{2.3.2.4}, hence in the weakly saturated class generated by $ \LH_P $.
\end{proof}


\subsection{Stability of weak equivalences under filtered colimits}\label{subsec:stabfilteredcolim}

In this subsection we explain how to fibrantly replace simplicial sets over $ P $ whose strata are Kan complexes without changing their strata, and use this to deduce that Joyal--Kan equivalences between such objects are Joyal equivalences (\Cref{prop:WEfiberwiseKan}).
We leverage this to show that Joyal--Kan equivalences are stable under filtered colimits (\Cref{lem:JKstableunderfilteredcolim}), verifying the last criterion to show that the Joyal--Kan model structure is simplicial (\Cref{thm:JKissimplicial}).
We deduce that the Joyal--Kan model structure presents the \category $ \StrP $ (\Cref{cor:underlying}).

\begin{notation}
	Let $ P $ be a poset.
	Write $ \Inv_P \subset \IH_P $ for those inner horn inclusions $ \incto{\horn{k}{n}}{\Delta^n} $ over $ P $ that are \textit{not vertical} in the sense that the stratification $ \fromto{\Delta^n}{P} $ is not a constant map.
\end{notation}

\begin{lemma}\label{lem:fibreplaceKan}
	Let $ f \colon \fromto{X}{Y} $ be a morphism in $ \sSetP $.
	Then there exists a commutative square
	\begin{equation*}
      \begin{tikzcd}
	       X \arrow[d, hooked, "i"'] \arrow[r, "f"] & Y \arrow[d, hooked, "j"] \\ 
	       \Xtilde \arrow[r, "\ftilde"'] & \Ytilde \comma
      \end{tikzcd}
    \end{equation*}
	in $ \sSetP $ where:
	\begin{enumerate}[label=\stlabel{lem:fibreplaceKan}, ref=\arabic*]
		\item\label{lem:fibreplaceKan.1} The morphisms $ i $ and $ j $ are $ \Inv_P $-cell maps.

		\item\label{lem:fibreplaceKan.2} The morphism $ \ftilde $ satisfies the right lifting property with respect to $ \Inv_P $.

		\item\label{lem:fibreplaceKan.3} The morphisms $ i $ and $ j $ restrict to isomorphisms on strata, i.e., for all $ p \in P $ the morphisms $ i $ and $ j $ restrict to isomorphisms of simplicial sets $ i \colon \isomto{X_p}{\Xtilde_p} $ and $ j \colon \isomto{Y_p}{\Ytilde_p} $.

		\item\label{lem:fibreplaceKan.4} If, in addition, all of the strata of $ X $ and $ Y $ are quasicategories, then $ \Xtilde $ and $ \Ytilde $ can be chosen to be quasicategories.
	\end{enumerate}

	In particular, if all of the strata of $ X $ and $ Y $ are Kan complexes, then $ \Xtilde $ and $ \Ytilde $ can be chosen to be fibrant in the Joyal--Kan model structure on $ \sSetP $.
\end{lemma}

\begin{proof}
	Since the morphisms in $ \Inv_P $ all have small domains, we can apply the small object argument to construct a square
	\begin{equation}\label{sq:Joyalreplacement}
      \begin{tikzcd}
	       X \arrow[d, hooked, "i"'] \arrow[r, "f"] & Y \arrow[d, hooked, "j"] \\ 
	       \Xtilde \arrow[r, "\ftilde"'] & \Ytilde \comma
      \end{tikzcd}
    \end{equation} 
    where $ i $ and $ j $ are $ \Inv_P $-cell maps and $ \ftilde $ has the right lifting property with respect to $ \Inv_P $, which proves \enumref{lem:fibreplaceKan}{1} and \enumref{lem:fibreplaceKan}{2}.

    To prove \enumref{lem:fibreplaceKan}{3} we examine the constructions of $ \Xtilde $ and $ \Ytilde $ via the small object argument.
    Both morphisms $ i $ and $ j $ are obtained by a transfinite composite of pushouts of inner horn inclusions $ \incto{\horn{k}{n}}{\Delta^n} $ in $ \Inv_P $.
    Hence to prove \enumref{lem:fibreplaceKan}{3} it suffices to show that given an object $ A \in \sSetP $ and a morphism $ f \colon \fromto{\horn{k}{n}}{A} $, where $ \horn{k}{n} \in \sSetP $ is the domain of a morphism $ g \colon \incto{\horn{k}{n}}{\Delta^n} $ in $ \Inv_P $, the morphism $ \gbar $ in the pushout square 
    \begin{equation*}
      \begin{tikzcd}
      	\horn{k}{n}  \arrow[dr, phantom, very near end, "\ulcorner", xshift=0.25em, yshift=-0.25em] \arrow[d, hooked, "g"'] \arrow[r, "f"] & A \arrow[d, hooked, "\gbar"] \\ 
       \Delta^n \arrow[r, "\fbar"'] & A'
      \end{tikzcd}
    \end{equation*}
    induces an isomorphism (of simplicial sets) on strata.
    To see this, let $ \sigma \colon \fromto{\Delta^n}{P} $ denote the stratification of the target of $ g $.
    Since $ g \in \Inv_P $, the stratification $ \sigma $ is not a constant map.
    We claim that for all $ p \in P $ and $ m \geq 0 $, the $ m $-simplices of $ A_p \subset A'_p $ and $ A'_p $ coincide.
    If $ p \notin \sigma(\Delta^n) $ or $ m < n-1 $, this is clear.
    Let us consider the remaining cases.
    \begin{enumerate}
    	\item If $ m = n - 1 $, then note that the only additional $ (n-1) $-simplex adjoined to $ A $ in the pushout defining $ A' $ is the image of the face
    	\begin{equation*}
    		\Delta^{\{0<\cdots<k-1<k+1<\cdots<n\}} \subset \Delta^n \period
    	\end{equation*}
    	Since the horn $ \horn{k}{n} \subset \Delta^n $ is an inner horn, both vertices $ \Delta^{\{0\}} $ and $ \Delta^{\{n\}} $ are contained in the face \smash{$\Delta^{\{0<\cdots<k-1<k+1<\cdots<n\}} $}.
    	Since the stratification $ \sigma \colon \fromto{\Delta^n}{P} $ is not constant, the image of the face \smash{$ \Delta^{\{0<\cdots<k-1<k+1<\cdots<n\}} $} in $ A' $ intersects more than one stratum.
    	Hence for each $ p \in P $, the $ (n-1) $-simplices of the strata $ A_p $ and $ A'_p $ coincide. 

    	\item If $ m = n $, then note that the only additional nondegenerate $ n $-simplex adjoined to $ A $ in the pushout defining $ A' $ is the unique nondegenerate $ n $-simplex of $ \Delta^n $.
    	Since the stratification $ \sigma \colon \fromto{\Delta^n}{P} $ is non-constant, the image of this top-dimensional simplex under $ \fbar $ intersects more than one stratum.
    	Similarly, note that since the image of the face \smash{$ \Delta^{\{0<\cdots<k-1<k+1<\cdots<n\}} $} in $ A' $ intersects more than one stratum (by the previous point), all of its degeneracies intersect more than one stratum. 
    	But the image of $ \Delta^n $ and images of the degeneracies of \smash{$ \Delta^{\{0<\cdots<k-1<k+1<\cdots<n\}} $} under $ \fbar $ are the only $ n $-simplices adoined to $ A $ in the pushout defining $ A' $.
    	Hence for each $ p \in P $, the $ n $-simplices of the strata $ A_p $ and $ A'_p $ coincide.

    	\item If $ m > n $, then the claim follows from the fact that the $ \el $-simplices of $ A_p $ and $ A'_p $ coincide for all $ \el \leq n $ and the $ n $-skeletality of $ \Delta^n $. 
    \end{enumerate}

    Now we prove \enumref{lem:fibreplaceKan}{4}; assume that the strata of $ X $ and $ Y $ are quasicategories. 
    To see that $ \Ytilde $ is a quasicategory, note that by the construction of the factorization \eqref{sq:Joyalreplacement} via the small object argument, $ \Ytilde $ is given by factoring the unique morphism $ \fromto{Y}{P} $ to the final object as a composite
    \begin{equation*}
      \begin{tikzcd}[sep=1.5em]
	       Y \arrow[r, hooked, "j"] & \Ytilde \arrow[r, "h"] & P 
      \end{tikzcd}
    \end{equation*} 
    of the $ \Inv_P $-cell map $ j $ followed by a morphism $ h $ with the right lifting property with respect to $ \Inv_P $.
    To show that $ \Ytilde $ is a quasicategory, we prove that $ h $ is an inner fibration.
    By the definition of $ \Inv_P $, the morphism $ h $ lifts against all inner horns \smash{$ \incto{\horn{k}{n}}{\Delta^n} $} over $ P $ where the stratification of $ \Delta^n $ is not constant.
    Thus to check that \smash{$ h \colon \fromto{\Ytilde}{P} $} is an inner fibration, it suffices to check that for every inner horn $ \incto{\horn{k}{n}}{\Delta^n} $ over $ P $ where the stratification $ \sigma \colon \fromto{\Delta^n}{P} $ is constant at a vertex $ p \in P $, every lifting problem 
    \begin{equation*}
      \begin{tikzcd}
	       \horn{k}{n} \arrow[d, hooked] \arrow[r] & \Ytilde \arrow[d, "h"] \\ 
	       \Delta^n \arrow[r, "\sigma"'] \arrow[ur, dotted] & P
      \end{tikzcd}
    \end{equation*}   
    admits a solution.
    The desired lift exists by \enumref{lem:fibreplaceKan}{3} because the stratum $ \Ytilde_p \isomorphic Y_p $ is a quasicategory by assumption.

    We conclude that $ \Xtilde $ is a quasicategory by showing that the stratification $ \fromto{\Xtilde}{P} $ is an inner fibration.
    First, note that since $ \ftilde \colon \fromto{\Xtilde}{\Ytilde} $ and $ h \colon \fromto{\Ytilde}{P} $ have the right lifting property with respect to $ \Inv_P $, so does the stratification $ h\ftilde \colon \fromto{\Xtilde}{P} $.
    Hence to show that the stratification $ \fromto{\Xtilde}{P} $ is an inner fibration, it suffices to show that $ \fromto{\Xtilde}{P} $ lifts against inner horns $ \incto{\horn{k}{n}}{\Delta^n} $ over $ P $ where the stratification $ \sigma \colon \fromto{\Delta^n}{P} $ is constant.
    Again, the desired lift exists by \enumref{lem:fibreplaceKan}{3} because the strata of $ \Xtilde $ are quasicategories.
\end{proof}


\begin{proposition}\label{prop:WEfiberwiseKan}
	Let $ P $ be a poset and $ f \colon \fromto{X}{Y} $ a morphism in \smash{$ \sSetP $}.
	If the strata of $ X $ and $ Y $ are all Kan complexes, then $ f $ is a Joyal--Kan equivalence if and only if $ f $ is a Joyal equivalence.
\end{proposition}

\begin{proof}
	Since the Joyal--Kan equivalences between fibrant objects of the Joyal--Kan model structure on $ \sSetP $ are precisely the Joyal equivalences \enumref{prop:JKexists}{7}, by 2-of-3 it suffices to show that there exists a commutative square
	\begin{equation*}
      \begin{tikzcd}
	       X \arrow[d, "i"', "\wr"{xshift=-0.1em}] \arrow[r, "f"] & Y \arrow[d, "j", "\wr"'{xshift=0.1em}] \\ 
	       \Xtilde \arrow[r, "\ftilde"'] & \Ytilde \comma
      \end{tikzcd}
    \end{equation*}
	in $ \sSetP $, where $ \Xtilde $ and $ \Ytilde $ are Joyal--Kan fibrant objects, and $ i \colon \equivto{X}{\Xtilde} $ and $ j \colon \equivto{Y}{\Ytilde} $ are Joyal equivalences.
	This follows from \Cref{lem:fibreplaceKan} since $ \Inv_P $-cell maps are, in particular, Joyal equivalences. 
\end{proof}

\begin{nul}\label{nul:filteredcolimobs}
	Since Joyal equivalences are stable under filtered colimits \cite[\HTTthm{Theorem}{2.2.5.1} \& \HTTpage{90}]{HTT}, to show that Joyal--Kan equivalences are stable under filtered colimits, \Cref{prop:WEfiberwiseKan} reduces us to constructing a functor \smash{$ F \colon \fromto{\sSetP}{\sSetP} $} that lands in strata-wise Kan complexes, admits a natural weak equivalence \smash{$ \equivto{\id{\sSetP}}{F} $}, and preserves filtered colimits.
	We accomplish this by applying Kan's $ \Ex^{\infty} $ functor vertically to each stratum.
\end{nul}

\begin{construction}\label{cnstr:badExP}
	Let $ P $ be a poset.
	Define a functor $ \ExPhat^{\infty} \colon \fromto{\sSetP}{\sSetP} $ by the assignment
	\begin{equation*}
		\goesto{X}{X \coproductlimits^{\Obj(P) \cross_P X} \Ex^{\infty}(\Obj(P) \cross_P X) \isomorphic X \coproductlimits^{\coprod_{p \in P} X_p} \paren{\textstyle \coprod_{p \in P} \Ex^{\infty}(X_p)}} \comma
	\end{equation*}
	where the pushout is taken in $ \sSetP $, and the stratifications
	\begin{equation*}
		\fromto{\coprod_{p \in P} \Ex^{\infty}(X_p)}{P} \andeq \fromto{\coprod_{p \in P} X_p}{P}
	\end{equation*}
	are induced by the constant maps $ \fromto{\Ex^{\infty}(X_p)}{P} $ and $ \fromto{X_p}{P} $ at $ p \in P $.
	
	We claim that the natural inclusion
	\begin{equation*}
		\incto{X}{\ExPhat^{\infty}(X)}
	\end{equation*}
	is a trivial Joyal--Kan cofibration.
	To see this, note that for each $ p \in P $, the inclusion $ \incto{X_p}{\Ex^{\infty}(X_p)} $ is a trivial Kan cofibration, so in the weakly saturated class generated by the horn inclusions $ \incto{\horn{k}{n}}{\Delta^n} $ in $ \sSet $, where $ n \geq 0 $ and $ 0 \leq k \leq n $.
	Thus, stratifying $ X_p $ and $ \Ex^{\infty}(X_p) $ via the constant maps at $ p \in P $, by \Cref{prop:stratifiedhornsall} the inclusion $ \incto{X_p}{\Ex^{\infty}(X_p)}  $ is a trivial Joyal--Kan cofibration in $ \sSetP $.
	To conclude, note that by definition the inclusion $ \incto{X}{\ExPhat^{\infty}(X)} $ is a pushout of the trivial Joyal--Kan cofibration
	\begin{equation*}
		\incto{\coprod_{p \in P} X_p}{\coprod_{p \in P} \Ex^{\infty}(X_p)} \rlap{\ .}
	\end{equation*}

	In particular, note that by \Cref{prop:WEfiberwiseKan} a morphism $ f \colon \fromto{X}{Y} $ in $ \sSetP $ is a Joyal--Kan equivalence if and only if 
	\begin{equation*}
		\ExPhat^{\infty}(f) \colon \fromto{\ExPhat^{\infty}(X)}{\ExPhat^{\infty}(Y)} 
	\end{equation*}
	is a Joyal equivalence.
\end{construction}

\begin{warning}
	The functor $ \ExPhat^{\infty} $ in general does \textit{not} preserve quasicategories over $ P $.
	In particular, if $ P $ is not discrete, then $ \ExPhat^{\infty} $ is not a fibrant replacement for the Joyal--Kan model structure.
\end{warning}

\begin{lemma}\label{lem:badExfilteredcolim}
	Let $ P $ be a poset.
	Then the functor $ \ExPhat^{\infty} \colon \fromto{\sSetP}{\sSetP} $ preserves filtered colimits.
\end{lemma}

\begin{proof}
	First, since filtered colimits commute with finite limits in $ \sSet $, the functor
	\begin{equation*}
		\Obj(P) \cross_P - \colon \fromto{\sSetP}{\sSetP}
	\end{equation*}
	preserves filtered colimits.
	Second, since Kan's $ \Ex^{\infty} $ functor preserves filtered colimits, the functor $ \fromto{\sSetP}{\sSetP} $ given by the assignment
	\begin{equation*}
		\goesto{X}{\Ex^{\infty}(\Obj(P) \cross_P X) \isomorphic \coprod_{p \in P} \Ex^{\infty}(X_p)}
	\end{equation*}
	preserves filtered colimits.
	The claim is now clear from the definition of $ \ExPhat^{\infty} $.
\end{proof}

Combining our observation \cref{nul:filteredcolimobs} with \Cref{lem:badExfilteredcolim} we deduce:

\begin{proposition}\label{lem:JKstableunderfilteredcolim}
	For any poset $ P $, Joyal--Kan equivalences in $ \sSetP $ are stable under filtered colimits.
\end{proposition}

\noindent \Cref{lem:JKstableunderfilteredcolim,2.rmk:pushout-prod,prop:incweakequiv} verify the three conditions of \Cref{nul:wantsimplcial} proving: 

\begin{theorem}\label{thm:JKissimplicial}
	For any poset $ P $, the Joyal--Kan model structure on $ \sSetP $ is simplicial.
\end{theorem}



\noindent From this we immediately deduce that the Joyal--Kan model structure presents $ \StrP $.

\begin{corollary}\label{cor:underlying}
	Let $ P $ be a poset.
	Then the underlying quasicategory of the Joyal--Kan model structure on $ \sSetP $ is the quasicategory $ \StrP $ of abstract $ P $-stratified homotopy types.
\end{corollary}


\section{The stratified homotopy hypothesis}\label{sec:JoyalKanandtopology}

The purpose of this section is to use our work on the Joyal--Kan model structure to prove our stratified homotopy hypothesis (\Cref{introprop:SingPequiv}).
In \cref{subsec:JoyalKanandtopology:elementaryresults}, we begin with some immediate consequences of our work in \cref{sec:modelonsSetP}.
\Cref{sec:strathomotopyhyp} proves the stratified homotopy hypothesis.
\Cref{subsec:global_SHH} explains a globalization where we allow the stratifying poset to vary.
In \cref{sec:otherwork}, we show that the Ayala--Francis--Tanaka--Rozenblyum homotopy theory of stratfied spaces embeds into ours.


\subsection{Elementary results}\label{subsec:JoyalKanandtopology:elementaryresults}

\begin{corollary}
	Let $ P $ be a poset and $ T $ a $ P $-stratified topological space.
	Then the $ P $-stratified simplicial set \smash{$ \SingP(T) $} is a Joyal--Kan fibrant object of \smash{$ \sSetP $} if and only if \smash{$ \SingP(T) $} is a quasicategory.
\end{corollary}

\begin{proof}
	Combine \Cref{nul:stratumofSingP,prop:equivfibrant}.
\end{proof}

Since the exit-path simplicial set of a conically stratified topological space is a quasicategory \Cref{nul:conicstratisqcat}, we deduce:

\begin{corollary}\label{cor:singPconicisfib}
	Let $ P $ be a poset.
	If \smash{$ T \in \TopP $} is conically stratified, then the simplicial set \smash{$ \SingP(T) $} is a Joyal--Kan fibrant object of $ \sSetP $.
\end{corollary}

\noindent Note that not all stratified topological spaces arising as geometric realizations of quasicategories are conically stratified:

\begin{example}\label{exm:notcon}
	Stratify the quasicategory $ \horn{0}{2} $ over $ [1] $ via the map sending $ 0 $ and $ 1 $ to $ 0 $  and $ 2 $ to $ 1 $.
	The $ [1] $-stratified topological space \smash{$ |\horn{0}{2}|_{[1]} $} is not conically stratified.
	Moreover, the $ [1] $-stratified simplicial set \smash{$ \Sing_{[1]} |\horn{0}{2}|_{[1]} $} is not a quasicategory.  
\end{example}

\begin{warning}
	\Cref{exm:notcon} shows that, unlike the Kan model structure, if $ P $ is a non-discrete poset, the functor $ \SingP \realP{-} $ is \textit{not} a fibrant replacement for the Joyal--Kan model structure on $ \sSetP $.
\end{warning}

\begin{notation}
	Let $ P $ be a poset.
	Write 
	\begin{equation*}
		\TopPex \subset \TopP
	\end{equation*}
	for the full subcategory spanned by those $ P $-stratified topological spaces $ T $ for which the simplicial set $ \SingP(T) $ is a quasicategory.
	Note that $ \TopPex $ contains all conically $ P $-stratified topological spaces.
\end{notation}

Since the Joyal--Kan model structure on $ \sSetP $ presents the \category $ \StrP $ (\Cref{cor:underlying}), we deduce the following variants of a result of Miller \cites[Theorem 6.10]{MillerThesis}[Theorem 6.3]{MR3066775}.

\begin{corollary}\label{cor:TopWEstrataandlink}
	Let $ P $ be a poset and let $ f \colon \fromto{T}{S} $ be a morphism in $ \TopPex $.
	Then the morphism \smash{$ \SingP(f) $} is an equivalence when regarded as a morphism in \smash{$ \StrP $} if and only if the following conditions are satisfied: 
	\begin{enumerate}[label=\stlabel{cor:TopWEstrataandlink}, ref=\arabic*]
		\item\label{cor:TopWEstrataandlink.1} For each $ p \in P $, the induced map on strata $ \fromto{T_p}{S_p} $ is a weak homotopy equivalence of topological spaces.

		\item\label{cor:TopWEstrataandlink.2} For all $ p,q \in P $ with $ p < q $, the induced map on topological links
		\begin{equation*}
			\fromto{\Map_{\TopP}(|\{p < q\}|_P,T)}{\Map_{\TopP}(|\{p < q\}|_P,S)}
		\end{equation*}
		is a weak homotopy equivalence of topological spaces.
	\end{enumerate}
\end{corollary}

\begin{proof}
	Combine \Cref{nul:stratumofSingP}, \Cref{nul:equivcheckedonstratandlink}, and \Cref{rem:Singoftoplink}.
\end{proof}

\begin{lemma}\label{lem:equivalencesforgoodspaces2}
	Let $ f \colon \fromto{T}{S} $ be a morphism in $ \TopPex $.
	The following conditions are equivalent:
	\begin{enumerate}[{\upshape (\ref*{lem:equivalencesforgoodspaces2}.1)}]
		\item The map $ f $ is a Douteau--Henriques equivalence.

		\item The morphism $ \SingP(f) $ is a weak equivalence in the Joyal--Kan model structure on $ \sSetP $.

		\item The morphism $ \SingP(f) $ is an equivalence when regarded as a morphism in the \category $ \StrP $ of abstract $ P $-stratified homotopy types. 
	\end{enumerate}
\end{lemma}

\begin{proof}
	Combine \Cref{nul:equivalencesforgoodspaces1,cor:underlying}.
\end{proof}


\subsection{The stratified homotopy hypothesis}\label{sec:strathomotopyhyp}

In this subsection we prove our stratified homotopy hypothesis (see \Cref{prop:SingPequiv}).
We accomplish this by appealing to the relationship between the \category $ \StrP $ and décollages explained in \cref{subsec:decollage,sec:accessible_loc}.

\begin{nul}
	Note that for any \smash{$ T \in \TopPex $}, the simplicial presheaf 
	\begin{equation*}
		\DP(T) \isomorphic \NerveP \SingP(T)
	\end{equation*}
	on $ \sd(P) $ introduced in \Cref{ntn:DP} already satisfies the Segal condition for décollages.
\end{nul}
	
\begin{notation}\label{ntn:WP}
	Let $ \WP \subset \Mor(\TopPex) $ denote the class of morphisms that induce weak homotopy equivalences on all strata and topological links.
	Equivalently, $ \WP $ is the class of morphisms that are sent to equivalences in the Joyal--Kan model structure under \smash{$ \SingP $} (\Cref{cor:TopWEstrataandlink,lem:equivalencesforgoodspaces2}).
\end{notation}	

\begin{nul}\label{comp:NLW}
	\Cref{thm:TopPpresentsStrP} implies that the induced functor of \categories
	\begin{align*}
		\ExitP \colon \TopPex[\WP^{-1}] &\to \StrP \\
		T &\mapsto \SingP(T)
	\end{align*}
	is fully faithful. 
	In fact $ \ExitP $ is an \textit{equivalence}:
\end{nul}

\begin{theorem}[(stratified homotopy hypothesis, local version)]\label{prop:SingPequiv}
	For any poset $ P $, the functor
	\begin{equation*}
		\ExitP \colon \fromto{\TopPex[\WP^{-1}]}{\StrP}
	\end{equation*}
	is an equivalence of \categories.
\end{theorem}

\begin{proof}
	By \Cref{comp:NLW} it suffices to show that $ \ExitP $ is essentially surjective.
	For this, it suffices to show that the induced functor
	\begin{equation*}
		\DP \colon \fromto{\TopPex[\WP^{-1}]}{\DecP}
	\end{equation*}
	is essentially surjective.
	We prove this by factoring the equivalence from a localization of the underlying \category of $ \TopP $ in the Douteau--Henriques model structure to $ \DecP $ through complete Segal spaces with a conservative functor to $ \Nerve(P) $ (cf. the proof of \Cref{thm:nerveequiv}).

	Let $ i \colon \incto{\sd(P)}{\DDelta_{/P}} $ denote the inclusion of the subdivision of $ P $ into the category of simplicies of $ P $.
	Write 
	\begin{equation*}
		\iupperstar \colon \fromto{\Fun(\Deltaop_{/P},\sSet)}{\Fun(\sdop{P},\sSet)}
	\end{equation*}
	for the restriction functor.
	Write $ \ilowershriek $ for the left adjoint of $ \iupperstar $, given by left Kan extension along $ i $.
	Then the induced adjunction
	\begin{equation}\label{eq:DecCSSadjunction}
		\adjto{\ilowershriek}{\sSetPproj}{\sSetDeltaPproj}{\iupperstar} \comma
	\end{equation}
	is a simplicial Quillen adjunction for the projective model structures (with respect to the Kan model structure on $ \sSet $).
	The simplicial Quillen equivalence 
	\begin{equation*}
		\adjto{\realP{-} \of \LP}{\sSetPproj}{\TopP}{\DP}  
	\end{equation*}
	factors as a composite of simplicial adjunctions
	\begin{equation*}
		\begin{tikzcd}[sep=2em]
			\sSetPproj \arrow[r, shift left, "\ilowershriek"] & \sSetDeltaPproj \arrow[r, shift left] \arrow[l, shift left, "\iupperstar"] & \TopP \comma \arrow[l, shift left, "\DP'"] 
		\end{tikzcd}
	\end{equation*}
	where the right adjoint $ \DP' \colon \fromto{\TopP}{\Fun(\DDelta_{/P}^{\op},\sSet)} $ is given by
	\begin{equation*}
		T \mapsto [(\Delta^n \to P) \mapsto \Sing \Map_{\TopP}(\realP{\Delta^n},T)] \period
	\end{equation*}
	Moreover, $ \DP' $ preserves Douteau--Henriques weak equivalences and fibrant objects, and $ \iupperstar $ preserves weak equivalences.
	Write
	\begin{equation*}
		\CSS_{/\Nerve(P)}^{\cons} \subset \CSS_{/\Nerve(P)}
	\end{equation*}
	for the full subcategory of the \category of complete Segal spaces over $ \Nerve(P) $ spanned by those complete Segal spaces $\fromto{C}{\Nerve(P)}$ such that for any $ p\in P $, the complete Segal space $ C_p $ is an \groupoid.
	Since the projective model structures on \smash{$ \Fun(\sdop{P},\sSet) $} and \smash{$ \Fun(\DDelta_{/P}^{\op},\sSet) $} are left proper, appealing to \HTT{Proposition}{A.3.7.8} we see that there are left Bousfield localization of \smash{$ \sSetPproj $} and \smash{$ \sSetDeltaPproj $} so that the Quillen adjunction induced by the Quillen adjunction \eqref{eq:DecCSSadjunction} presents the equivalence $ \DecP \equivalent \CSS_{/\Nerve(P)}^{\cons} $ from the proof of \Cref{thm:nerveequiv}.

	Let $ \WP' \subset \Mor(\TopP) $ denote the class of morphisms sent by $ \DP' $ to weak equivalences in the left Bousfield localization of \smash{$ \sSetDeltaPproj $} presenting the \category $ \CSS_{/\Nerve(P)}^{\cons} $.
	From \Cref{thm:Douteaumain} we deduce that $ \DP' $ descends to an equivalence of \categories
	\begin{equation}\label{eq:DPequiv}
		\DP' \colon \equivto{\TopP[(\WP')^{-1}]}{\CSS_{/\Nerve(P)}^{\cons}} \period
	\end{equation}
	Thus the equivalence \eqref{eq:DPequiv} restricts to an equivalence
	\begin{equation*}
		\equivto{\TopPseg[(\WP')^{-1}]}{\CSS_{/\Nerve(P)}^{\cons}} \comma
	\end{equation*}
	where \smash{$ \TopPseg \subset \TopP $} is the full subcategory spanned by those objects $ T $ sent to complete Segal spaces under $ \DP' $.
	The functor 
	\begin{equation*}
		\DP' \colon \fromto{\TopP}{\Fun(\DDelta_{/P}^{\op},\sSet)}
	\end{equation*}
	is the composite of \smash{$ \SingP \colon \fromto{\TopP}{\sSetP} $} with the `nerve' functor
	\begin{align*}
		\Nerve'_P \colon \sSetP & \to \Fun(\DDelta_{/P}^{\op},\sSet) \\
		X &\mapsto [(\Delta^n \to P) \mapsto \Map_{\sSetP}(\Delta^n,X)] \period
	\end{align*}
	By \cite[Corollary 3.6]{MR2342834} we see that $ \DP'(T) $ is a complete Segal space if and only if $ \SingP(T) $ is a quasicategory.
	Hence \smash{$ \TopPseg = \TopPex $}.
	To conclude, note that the class of morphisms in $ \TopPex $ that lie in $ \WP' $ coincides with the class of morphisms sent to Joyal--Kan equivalences under $ \SingP $.
\end{proof}


\subsection{Changing the stratifying poset}\label{subsec:global_SHH}

In this subsection, we explain why our stratified homotopy hypothesis for a fixed poset (\Cref{prop:SingPequiv}) implies a `global' version where we allow the stratifying poset to vary.
We first fix some notation.

\begin{notation}
	Write
	\begin{equation*}
		\StrTop \subset \Fun([1],\Top)
	\end{equation*}
	for the full subcategory of the arrow category of $ \Top $ on those morphisms $ s \colon \fromto{T}{P} $ where $ P $ is (the Alexandroff topological space associated to) a poset, and the simplicial set $ \SingP(T) $ is a quasicategory.

	Similarly, write
	\begin{equation*}
		\Str \subset \Fun([1],\Catinfty)
	\end{equation*}
	for the full subcategory of the arrow category of $ \Catinfty $ on those functors $ f \colon \fromto{C}{P} $ where $ P $ is a poset, and $ f $ is conservative.

	Write $ \Exit $ for the functor
	\begin{align*}
		\StrTop &\to \Str \\
		[T \to P] &\mapsto [\ExitP(T) \to P] \period
	\end{align*}
	Write $ \W \subset \Mor(\StrTop) $ for the set of morphisms in $ \StrTop $ sent to equivalences by the functor $ \Exit $.
\end{notation}

\begin{nul}
	Note that a morphism 
	\begin{equation*}
		\begin{tikzcd}
			T \arrow[r, "f"] \arrow[d] & S \arrow[d] \\ 
			P \arrow[r, "\phi"'] & Q \comma
		\end{tikzcd}
	\end{equation*}
	is in $ \W $ if and only if $ \phi $ is an isomorphism and $ f $ induces weak homotopy equivalences on all strata and links (\Cref{cor:TopWEstrataandlink}).
\end{nul}

The following is an immediate consequence of \Cref{prop:SingPequiv}.

\begin{corollary}[(stratified homotopy hypothesis, global version)]\label{cor:globalSHH}
	The functor
	\begin{equation*}
		\Exit \colon \fromto{\StrTop[\W^{-1}]}{\Str}
	\end{equation*}
	is an equivalence of \categories.
\end{corollary}


\subsection{Relation to conically smooth stratified spaces}\label{sec:otherwork}

We conclude by using our stratified homotopy hypothesis to show that the Ayala--Francis--Tanaka--Rozenblyum homotopy theory of conically smooth stratfied spaces embeds into our homotopy theory of stratified spaces.

\begin{recollection}[(conically smooth stratified spaces)]
	In work with Tanaka \cite[\S3]{MR3590534}, Ayala and Francis introduced \textit{conically smooth structures} on stratified topological spaces.
	Ayala and Francis furtuer studied these in work with Rozenblyum \cite{MR3941460}.
	Write $ \Con $ for their category of conically smooth stratified spaces and conically smooth maps.
	An object of $ \Con $ consists of a stratified topological space $ s \colon \fromto{T}{P} $ where
	\begin{enumerate}[label=\stlabel{comp:AFR}]
		\item for each $ p \in P $, the set $ \setbar{q \in P}{q \leq p} $ is finite, and

		\item each stratum of $ s $ is connected,
	\end{enumerate}
	along with an additional \textit{structure} of a `conically smooth atlas' on $ T $.
	The conically smooth atlas, in particular, endows the strata of $ T $ with the structure of smooth manifolds.
	Morphisms in $ \Con $ are commutative squares of continuous maps
	\begin{equation*}
		\begin{tikzcd}
			T \arrow[r, "f"] \arrow[d] & S \arrow[d] \\ 
			P \arrow[r] & Q \comma
		\end{tikzcd}
	\end{equation*}
	where the map $ f $ satisfies additional regularity hypotheses.
	When $ P = Q = \ast $, these regularity hypotheses require that $ f $ be a smooth map of smooth manifolds.  

	Conically smooth stratified spaces are, in particular, conically stratified.
	Hence there is a functor
	\begin{equation*}
		\fromto{\Con}{\StrTop}
	\end{equation*}
	forgetting the conically smooth structure.
	The Ayala--Francis--Rozenblyum \category of stratified spaces is the \category obtained from $ \Con $ by inverting the class $ \Hup $ of stratified homotopy equivalences \cite[Theorem 2.4.5]{MR3941460}.
\end{recollection}

\begin{comparison}\label{comp:AFR}
	The composite functor
	\begin{equation*}
		\begin{tikzcd}
			\Con \arrow[r] & \StrTop \arrow[r, "\Exit"] & \Str
		\end{tikzcd}
	\end{equation*}
	sends the class $ \Hup $ to equivalences, hence induces a functor of \categories
	\begin{equation*}
		\Exit \colon \fromto{\Con[\Hup^{-1}]}{\Str} \period
	\end{equation*}
	As a result of \cite[Lemma 3.3.9 \& Theorem 4.2.8]{MR3941460} the functor \smash{$ \Exit \colon \fromto{\Con[\Hup^{-1}]}{\Str} $} is fully faithful.
	Thus we have a commutative triangle of fully faithful functors of \categories
	\begin{equation*}
		\begin{tikzcd}[column sep={20ex,between origins}, row sep={8ex,between origins}]
			\Con[\Hup^{-1}] \arrow[d, hooked] \arrow[r, hooked, "\Exit"] & \Str \\ 
			\StrTop[\W^{-1}] \arrow[ur, shift right=0.25em, "\Exit"', "\sim"{sloped, xshift=-0.1em, yshift=-0.15em}] & \phantom{\Str} \comma
		\end{tikzcd}
	\end{equation*}
	where the vertical functor is induced by the forgetful functor $ \fromto{\Con}{\StrTop} $.
	In particular, the theory of conically stratified spaces with equivalences on exit-path \categories inverted subsumes the Ayala--Francis--Tanaka--Rozenblyum theory of stratified spaces.
\end{comparison}

\begin{remark}
	Ayala--Francis--Rozenblyum conjectured that every Whitney-stratified space admits a conically smooth structure \cite[Conjecture 0.0.7]{MR3941460}.
	For a long time this conjecture was open; Nocera and Volpe have recently proven it \cite[Theorem 2.7]{arXiv:2105.09243}.
\end{remark}



\DeclareFieldFormat{labelnumberwidth}{#1}
\printbibliography[keyword=alph, heading=references]
\DeclareFieldFormat{labelnumberwidth}{{#1\adddot\midsentence}}
\printbibliography[heading=none, notkeyword=alph]

\end{document}